\documentclass[a4paper,12pt]{article}

\makeatletter
 
  \@addtoreset{equation}{section}
 \makeatother
\usepackage{amsfonts}
\usepackage{amssymb}
\usepackage{mathrsfs}
\usepackage{amsmath,amsthm}
\usepackage{amsmath}
\usepackage{algorithm}
\usepackage{algorithmic}
\usepackage{indentfirst}
\usepackage{float}
\usepackage{physics}
\usepackage[pdftex]{graphicx}

\begin{document}
\setlength{\textwidth}{6.5in}
\setlength{\oddsidemargin}{0in}
\setlength{\topmargin}{-0.52in}
\setlength{\textheight}{9.0in}
\setlength{\footskip}{0.7in}

\newtheorem{definition}{Definition}
\newtheorem{assumption}{$[$ A}
\newtheorem{condition}{$[$ C}
\newtheorem{lemma}{Lemma}
\newtheorem{proposition}{Proposition}
\newtheorem{theorem}{Theorem}
\newtheorem{remark}{Remark}
\newtheorem{example}{Example}
\def\n{{\bf n}}
\def\A{{\bf A}}
\def\B{{\bf B}}
\def\C{{\bf C}}
\def\D{{\bf D}}
\def\E{{\bf E}}
\def\e{{\bf e}}
\def\F{{\bf F}}
\def\G{{\bf G}}
\def\H{{\bf H}}
\def\I{{\bf I}}
\def\J{{\bf J}}
\def\K{{\bf K}}
\def\L{{\bf L}}
\def\M{{\bf M}}
\def\N{{\bf N}}
\def\O{{\bf O}}
\def\P{{\bf P}}
\def\Q{{\bf Q}}
\def\R{{\bf R}}
\def\S{{\bf S}}
\def\T{{\bf T}}
\def\U{{\bf U}}
\def\V{{\bf V}}
\def\W{{\bf W}}
\def\X{{\bf X}}
\def\Y{{\bf Y}}
\def\Z{{\bf Z}}
\def\cala{{\cal A}}
\def\calb{{\cal B}}
\def\calc{{\cal C}}
\def\cald{{\cal D}}
\def\cale{{\cal E}}
\def\calf{{\cal F}}
\def\calg{{\cal G}}
\def\calh{{\cal H}}
\def\cali{{\cal I}}
\def\calj{{\cal J}}
\def\calk{{\cal K}}
\def\call{{\cal L}}
\def\calm{{\cal M}}
\def\caln{{\cal N}}
\def\calo{{\cal O}}
\def\calp{{\cal P}}
\def\calq{{\cal Q}}
\def\calr{{\cal R}}
\def\cals{{\cal S}}
\def\calt{{\cal T}}
\def\calu{{\cal U}}
\def\calv{{\cal V}}
\def\calw{{\cal W}}
\def\calx{{\cal X}}
\def\caly{{\cal Y}}
\def\calz{{\cal Z}}
%
\def\sskip{\hspace{0.5cm}}
\def\simleq{ \raisebox{-.7ex}{\em $\stackrel{{\textstyle <}}{\sim}$} }
\def\leqsim{ \raisebox{-.7ex}{\em $\stackrel{{\textstyle <}}{\sim}$} }
\def\ep{\epsilon}
\def\vep{\varepsilon}
\def\half{\frac{1}{2}}
\def\iku{\rightarrow}
\def\Iku{\Rightarrow}
\def\ikup{\rightarrow^{p}}
\def\inclusion{\hookrightarrow}
\def\cadlag{c\`adl\`ag\ }
\def\up{\uparrow}
\def\down{\downarrow}
\def\doti{\Leftrightarrow}
\def\douti{\Leftrightarrow}
\def\dochi{\Leftrightarrow}
\def\douchi{\Leftrightarrow}%
\def\yy{\\ && \nonumber \\}
\def\y{\vspace*{3mm}\\}
\def\nn{\nonumber}
\def\be{\begin{equation}}
\def\ee{\end{equation}}
\def\bea{\begin{eqnarray}}
\def\eea{\end{eqnarray}}
\def\beas{\begin{eqnarray*}}
\def\eeas{\end{eqnarray*}}
%
\def\ha{\hat{A}}
\def\hc{\hat{c}}
\def\hd{\hat{D}}
\def\hv{\hat{V}}
\def\hsd{{\hat{d}}}
\def\hx{\hat{X}}
\def\hsx{\hat{x}}
\def\hxa{\widehat{XA}}

\def\bsx{\bar{x}}
\def\bsd{{\bar{d}}}
\def\bx{\bar{X}}
\def\ba{\bar{A}}
\def\bb{\bar{B}}
\def\bc{\bar{C}}
\def\bv{\bar{V}}

\def\balpha{\bar{\alpha}}
\def\bbalpha{\bar{\bar{\alpha}}}
\def\combi{\l(\begin{array}{c}\alpha\\ \beta \end{array}\r)}
\def\f{^{(1)}}
\def\s{^{(2)}}
\def\ss{^{(2)*}}
\def\l{\left}
\def\r{\right}
\def\a{\alpha}
\def\b{\beta}
\def\L{\Lambda}
\def\lam{\lambda}

\newtheorem{thm}{Theorem}
\newtheorem{prop}{Proposition}
\newtheorem{defn}{Definition}
\newtheorem{rem}{Remark}
\newtheorem{step}{Step}
\newtheorem{cor}{Corollary}
\newtheorem{assump}{Assumption}

\everymath {\displaystyle}

\newcommand{\ruby}[2]{
\leavevmode
\setbox0=\hbox{#1}
\setbox1=\hbox{\tiny #2}
\ifdim\wd0>\wd1 \dimen0=\wd0 \else \dimen0=\wd1 \fi
\hbox{
\kanjiskip=0pt plus 2fil
\xkanjiskip=0pt plus 2fil
\vbox{
\hbox to \dimen0{
\small \hfil#2\hfil}
\nointerlineskip
\hbox to \dimen0{\mathstrut\hfil#1\hfil}}}}

\def\qedsymbol{$\blacksquare$}
\renewcommand{\thefootnote }{\fnsymbol{footnote}}
\renewcommand{\refname }{References}

\everymath {\displaystyle}

\allowdisplaybreaks[1]

\title{
A new efficient approximation scheme for solving high-dimensional semilinear PDEs: 
control variate method for Deep BSDE solver 
}
\author{Akihiko Takahashi\footnote{The University of Tokyo, Tokyo, Japan}, \ Yoshifumi Tsuchida\footnote{Hitotsubashi University, Tokyo, Japan} \ and Toshihiro Yamada\footnote{Hitotsubashi University, Tokyo, Japan} \footnote{Japan Science and Technology Agency (JST), Tokyo, Japan}}
\date{First version: January 21, 2021, \  This version: January 30, 2021}
\maketitle

\abstract
This paper introduces a new approximation scheme for solving high-dimensional semilinear partial differential equations (PDEs) and backward stochastic differential equations (BSDEs). First, we decompose a target semilinear PDE (BSDE) into two parts, namely ``dominant" linear and ``small" nonlinear PDEs.
Then, we employ a Deep BSDE solver with a new control variate method to solve those PDEs, where approximations based on an asymptotic expansion technique are effectively applied to the linear part and also used as control variates for the nonlinear part.
Moreover, our theoretical result indicates that errors of the proposed method become much smaller than those of the original Deep BSDE solver. Finally, we show numerical experiments 
to demonstrate the validity of our method, which is consistent with the theoretical result in this paper.

\section{Introduction}
High-dimensional semilinear partial differential equations (PDEs) are often used to describe various complex, large-scale phenomena appearing in physics, applied mathematics, economics and finance. Such PDEs typically have the form: 
\begin{eqnarray}
	& \frac{\partial}{\partial t}u(t,x)+{\cal L}u(t,x)+f(t,x,u(t,x),\partial_xu(t,x)\sigma(t,x))=0, \ \ \ t<T, \ \ x\in \mathbb{R}^d,   
	\label{semilinearPDE}\\
	& u(T,x)=g(x),  \ \ x\in \mathbb{R}^d, \nn
\end{eqnarray} 
where $f$ is a nonlinear function, ${\cal L}$ is a second order differential operator of the type: 
\begin{eqnarray}
	{\cal L}\varphi(t,x)=\sum_i \mu^i(t,x) \partial_{x_i}\varphi(t,x)+\frac{1}{2}\sum_{i,j}[\sigma\sigma^{\top}]_{i,j}(t,x)\partial_{x_i}\partial_{x_j} \varphi(t,x), 
\end{eqnarray} 
and the dimension $d$ is assumed to be high. To solve the nonlinear PDE, we have to rely on some numerical schemes since they have no closed-form solutions especially in high-dimensional cases. Classical methods such as finite differences and finite elements fail in high-dimensional cases due to their exponential growth of complexity. 
In the last two decades, probabilistic approaches have been studied with Monte Carlo methods for backward stochastic differential equations (BSDEs) since  solutions of semilinear PDEs can be represented by the ones of corresponding BSDEs through the nonlinear Feynman-Kac formula (see Zhang (2017) \cite{Z} for instance).  

In Weinan E et al. (2017) \cite{EHJ}, a novel computational scheme called the {\it Deep BSDE method} is proposed. In the Deep BSDE method, a stochastic target problem is considered with a forward-discretization scheme of the related BSDE. Then, the control problem is solved with a deep learning algorithm. The Deep BSDE method has opened the door to tractability of higher dimensional problems,
which enables us to solve high-dimensional semilinear PDEs within realistic computation time. 
Recently,  notable related works, mostly with neural networks have developed new methods for solving various types of high dimensional PDEs. See \cite{BGJ}\cite{EHJ20}\cite{EGJS}\cite{GRP}\cite{GJS}\cite{GHJZ}\cite{HL}\cite{HLZ}\cite{HJS}\cite{HJKNW}\cite{HZW}\cite{LLM}\cite{SS}\cite{ZBYZ} for example.  

While high-dimensional semilinear PDEs can be feasibly solved by the Deep BSDE method, the deviation 
of its estimated value from the true one is not small with reasonable computational time. Then, constructing an acceleration scheme for the Deep BSDE method is desirable.

Fujii et al. (2019) \cite{FTT} proposed an improved scheme for the Deep BSDE method. They used a prior knowledge with an asymptotic expansion method for a target BSDE and obtained its fast approximation. Then, they found that numerical errors become small in accordance with the fast decrease in values of the corresponding loss function. The scheme enables us to reduce processing load of the original Deep BSDE solver. For details of the asymptotic expansion method, a key technique applied in their article,
see Takahashi (1999, 2015) \cite{T1}\cite{T2}, Kunitomo and Takahashi (2001, 2003) \cite{KT1}\cite{KT2} and references therein. Moreover, Naito and Yamada (2020) \cite{NY} presented an extended scheme of Fujii et al. (2019) \cite{FT} by applying the backward Euler scheme for a BSDE with a good initial detection of the solution to a target PDE so that the Deep BSDE method works more efficiently. 

In the current work, we develop a new deep learning-based approximation for solving high-dimensional semilinear PDEs by extending the schemes in Weinan E et al. (2017) \cite{EHJ}, Fujii et al. (2019) \cite{FTT} and Naito and Yamada (2020) \cite{NY}. In particular, we propose an efficient control variate method for the Deep BSDE solver in order to obtain more accurate and stable approximations. Let us briefly explain the strategy considered in this paper. We first decompose the semilinear PDE (\ref{semilinearPDE}) into two parts,

\noindent
$u(t,x)= {\cal U}^1(t,x)+ {\cal U}^2(t,x)$ as follows: 
\begin{eqnarray}
	& \frac{\partial}{\partial t}{\cal U}^1(t,x)+{\cal L}{\cal U}^1(t,x)=0, \ \ \ t<T, \ \ x\in \mathbb{R}^d,  \label{LPDE}\\
	& {\cal U}^1(T,x)=g(x), \ \  x\in \mathbb{R}^d, \nn
\end{eqnarray} 
and 
\begin{align}
	& \frac{\partial}{\partial t}{\cal U}^2(t,x)+{\cal L}{\cal U}^2(t,x) \label{NLPDE} \\
	& \hspace{0.5em} +f(t,x,{\cal U}^1(t,x)+{\cal U}^2(t,x),\partial_x{\cal U}^1(t,x)\sigma(t,x)+\partial_x{\cal U}^2(t,x)\sigma(t,x))=0, \ \  t<T, \  x\in \mathbb{R}^d, \nn \\
	& \hspace{10.0em} {\cal U}^2(T,x)=0,  \ \ x\in \mathbb{R}^d. \nn
\end{align} 
Here, we remark that the solution $u$ of the semilinear PDE (\ref{semilinearPDE}) is given by the sum of the solutions ${\cal U}^1$ and ${\cal U}^2$ of PDEs (\ref{LPDE}) and (\ref{NLPDE}), respectively. Also, we note that ${\cal U}^1$ is the solution to the ``dominant" linear PDE and ${\cal U}^2$ is the solution to the ``small" 
residual nonlinear PDE with null terminal condition whose magnitude is governed by the driver $(t,x,y,z) \mapsto f(t,x,{\cal U}^1(t,x)+y,\partial_x{\cal U}^1(t,x)\sigma(t,x)+z)$, which is 
generally expected 
to have small nonlinear effects on the solution of $u$. 
Consequently, 
the decomposition of the target $u(0,\cdot)$ is represented as follows: 
\begin{align}
	u(0,x)=\underset{\mbox{\tiny{``dominant" linear PDE part}}}{{\cal U}^1(0,x)}+\underset{\mbox{\tiny{``small" nonlinear PDE part}}}{{\cal U}^2(0,x)}, \ \ \ \ x\in \mathbb{R}^d. \label{decomposition}
\end{align} 

We next approximate 
\begin{enumerate}
	\item ${\cal U}^1$ by an asymptotic expansion method denoted by
	${\cal U}^{1,\mathrm{Asymp}}$;
	\item ${\cal U}^2$ by the Deep BSDE method, denoted by
	${\cal U}^{2,\mathrm{Deep}}$.
\end{enumerate}
We expect that ${\cal U}^{1,\mathrm{Asymp}}$ in the approximation
\begin{align}
	u(0,x)\approx {\cal U}^{1,\mathrm{Asymp}}(0,x)+{\cal U}^{2,\mathrm{Deep}}(0,x),\ \ \ \ x\in \mathbb{R}^d, \label{decomp-appx}
\end{align} 
becomes a control variate. Furthermore, ${\cal U}^{1,\mathrm{Asymp}}$ and $\partial_x{\cal U}^{1,\mathrm{Asymp}}\sigma$ in the approximate driver $(t,x,y,z) \mapsto f(t,x,{\cal U}^{1,\mathrm{Asymp}}(t,x)+y,\partial_x{\cal U}^{1,\mathrm{Asymp}}(t,x)\sigma(t,x)+z)$ of ${\cal U}^{2,\mathrm{Deep}}$ will be doubly the control variates. The current work shows how the proposed method works well
as a new deep learning-based approximation in both theoretical and numerical aspects. 

The organization of this paper is as follows: The next section briefly introduces the deep BSDE solver and acceleration schemes with asymptotic expansions. Section 3 explains our proposed method with the main theoretical result and Section 4 presents our numerical scheme with its experiment.

\section{Deep BSDE solver and acceleration scheme with asymptotic expansion}
Let $T>0$ and $(\Omega,{\cal F}, \{ {\cal F}_t \}_{0\leq t \leq T},P)$ be a filtered probability space equipped with a $d$-dimensional Brownian motion $W=\{ (W_t^1,\cdots,W_t^d) \}_{0 \leq t\leq T}$ and a square-integrable $\mathbb{R}^d$-valued random variable $\xi$, which is independent of $W$. The filtration $\{ {\cal F}_t \}_{0\leq t \leq T}$ is generated by $\{W_t + \xi \}_{0\leq t\leq T}$. Under this setting we consider the following FBSDE: 
\begin{align}
  dX_t^\varepsilon
  =&\mu(t,X_t^\varepsilon)dt+\varepsilon \sigma(t,X_t^\varepsilon)dW_t,\quad
  X_0^\varepsilon=\xi,\label{FSDE}\\
  -dY_t^{\varepsilon,\alpha}
  =&\alpha f(t, X_t^\varepsilon, Y_t^{\varepsilon,\alpha}, Z_t^{\varepsilon,\alpha}) dt - Z_t^{\varepsilon,\alpha} dW_t,\quad
  Y_T^{\varepsilon,\alpha}=g(X_T^\varepsilon), \label{BSDE}
\end{align}
where $\mu$ is a $\mathbb{R}^d$-valued function on $[0,T] \times \mathbb{R}^d$, $\sigma$ is a $\mathbb{R}^{d \otimes d}$-valued function on $[0,T] \times \mathbb{R}^d$, $f: [0,T] \times \mathbb{R}^d \times \mathbb{R} \times \mathbb{R}^d \to \mathbb{R}$, $g:\mathbb{R}^d \to \mathbb{R}$ are some functions so that the FBSDE has the unique solution, and $\varepsilon, \alpha \in (0,1)$ are some small parameters. 
Here, we assume that $\mu$ and $\sigma$ are bounded and smooth in $x$ and have bounded derivatives with any orders. Also, $f$ is uniformly Lipschitz continuous function with the Lipschitz constant $C_{\mathrm{Lip}}[f]$ and at most linear growth in the variables $x,y,z$. The function $g$ is assumed to be $C^2_b$-class. The functions $\mu,\sigma,f$ are uniformly H\"older-$1/2$ continuous with respect to $t$. Furthermore, we put the condition that there is $\varepsilon_0>0$ such that $\sigma(t,x)\sigma(t,x)^{\top} \geq \varepsilon_0 I$ for all $t \in [0,T]$ and $x \in \mathbb{R}^N$. We sometimes omit the subscripts ${\cdot}^{\varepsilon}$ or ${\cdot}^{\varepsilon,\alpha}$ 
if no confusion arises.  

 The corresponding semilinear PDE is given by
\begin{gather}
  \partial_t u(t,x) + \mathcal{L}^\varepsilon u(t,x) +  f^{\alpha}(t,x, u(t,x), \partial_x u(t,x) \sigma^{\varepsilon}(t,x))=0, \ \ \ t<T, \\
  u(T,x)=g(x),\nn
\end{gather}
where $\sigma^{\varepsilon}=\varepsilon \sigma$, $f^{\alpha}=\alpha f$, $\partial_x=(\partial_{x_{1}},\cdots,\partial_{x_{d}})=({\partial}/{\partial x_{1}},\cdots,{\partial}/{\partial x_{d}})$ and ${\cal L}^\varepsilon$ is the generator: 
\begin{eqnarray}
{\cal L}^\varepsilon =\sum_{i=1}^d\mu^i(t,x) \frac{\partial}{\partial x_i}+\frac{1}{2} \sum_{i_1,i_2=1}^d \sigma^{\varepsilon,i_1}(t,x)\sigma^{\varepsilon,i_2}(t,x) \frac{\partial^2}{\partial x_{i_1}\partial x_{i_2}}.
\end{eqnarray}
The purpose of this paper is to estimate 
\begin{eqnarray}
u(0,X_0^\varepsilon)=Y_0^{\varepsilon,\alpha}
\end{eqnarray}
against high dimensional FBSDEs/semilinear PDEs. In particular, we introduce an approximation with a deep BSDE solver to propose an efficient control variate method for solving semilinear PDEs. To explain how our method works well as a new scheme, we briefly review the deep BSDE method proposed in Weinan E et al. (2017) \cite{EHJ} and an approximation method developed by Fujii et al. (2019) \cite{FTT}. 

\subsection{Deep BSDE method by Weinan E et al. (2017)}
In Weinan E et al. (2017) \cite{EHJ}, the authors considered the minimization problem of the loss function: 
\begin{eqnarray}
\inf_{ Y_0^{\varepsilon,\alpha,(n)},  {Z}^{\varepsilon,\alpha,(n)}} 
\Big\| g(\bar{X}_T^{\varepsilon,(n)})-{Y}_T^{\varepsilon,\alpha,(n)} \Big\|_2^2 \label{original_DBSDE_optimization_prob}
\end{eqnarray}
where $\| \cdot \|_2=E[| \cdot |^2]^{1/2}$, subject to
\begin{align}
 {Y}_t^{\varepsilon,\alpha,(n)}&= {Y}^{\varepsilon,\alpha,(n)}_0 -\int_0^t  f^{\alpha}(s,\bar{X}_s^{\varepsilon,(n)},Y_s^{\varepsilon,\alpha,(n)},Z_s^{\varepsilon,\alpha,(n)})ds+\int_0^t  Z_s^{\varepsilon,\alpha,(n)} dW_s,
\end{align}
where $\bar{X}^{\varepsilon,(n)}$ is the continuous Euler-Maruyama scheme
with number of discretization time steps $n$: 
\begin{align}
\bar{X}_t^{\varepsilon,(n)}=\int_0^t \mu(\varphi(s),\bar{X}_{\varphi(s)}^{\varepsilon,(n)})ds+\int_0^t \sigma^{\varepsilon}(\varphi(s),\bar{X}_{\varphi(s)}^{\varepsilon,(n)})dW_s, \ \ \ t \geq 0, 
\end{align}
with $\varphi(s)=\max \{ kT/n; \ s \geq kT/n \}$. They solved the problem by using a deep learning algorithm and checked the effectiveness of the method for nonlinear BSDEs/PDEs even for the high dimension $d$. The method is known as Deep BSDE solver. 

Then, we have 
\begin{align}
 {Y}_0^{\varepsilon,\alpha} \ \approx \ {Y}^{\varepsilon,\alpha,(n),\ast}_0, 
\end{align}
where $Y_0^{\varepsilon,\alpha,(n),\ast}$ is obtained by solving (\ref{original_DBSDE_optimization_prob}), which is justified by the following estimate shown in Han and Long (2020)~\cite{HL}. 

\begin{thm}[Han and Long (2020)]\label{HL-thm1}
  There exists $C>0$ such that 
  \begin{equation}
      E[|Y_0^{\varepsilon,\alpha} - Y_0^{\varepsilon,\alpha,(n)}|^2]
      \leq C \frac{1}{n} + C \norm{g(\bar{X}_T^{\varepsilon,(n)}) - {Y}_T^{\varepsilon,\alpha,(n)}}_2^2,\label{han-long-thm}
  \end{equation}
  for $n\geq 1$. 
\end{thm}

\subsection{An approximation method by Fujii et al. (2019)}
In Fujii et al. (2019) \cite{FTT}, the authors considered the problem 
\begin{eqnarray}
\inf_{ \widetilde{Y}_0^{\varepsilon,\alpha,(n)},  \widetilde{Z}^{\mathrm{Res},\varepsilon,\alpha,(n)}} \Big\| g(\bar{X}_T^{\varepsilon,(n)})-\widetilde{Y}_T^{\varepsilon,\alpha,(n)} \Big\|_2^2 \label{FTT_DBSDE_optimization_prob}
\end{eqnarray}
subject to
\begin{align}
\widetilde{Y}_t^{\varepsilon,\alpha,(n)}&= \widetilde{Y}^{\varepsilon,\alpha,(n)}_0 -\int_0^t  f^{\alpha}(s,\bar{X}_s^{\varepsilon,(n)},\widetilde{Y}_s^{\varepsilon,\alpha,(n)},\widehat{Z}_s^{\varepsilon,\alpha,(n)}+\widetilde{Z}_s^{\mathrm{Res},\varepsilon,\alpha,(n)})ds\\
&\quad+\int_0^t  \{ \widehat{Z}_s^{\varepsilon,\alpha,(n)}+\widetilde{Z}_s^{\mathrm{Res},\varepsilon,\alpha,(n)} \} dW_s,
\end{align}
where $\widehat{Z}^{\varepsilon,\alpha,(n)}$ is a prior knowledge of $Z$ which is easily computed by an asymptotic expansion method, and they solve the minimization problem 
with respect to $\widetilde{Y}_0^{\varepsilon,\alpha,(n)}$ and $\widetilde{Z}^{\mathrm{Res},\varepsilon,\alpha,(n)}$ by Deep BSDE solver. The authors showed that the scheme gives better accuracy than the original Deep BSDE solver. Furthermore, Naito and Yamada (2020) \cite{NY} proposed an acceleration scheme by extending the method of Fujii et al. (2019) \cite{FTT} with a good initial detection of $Y_0$ and the backward Euler scheme of $Z$. They confirmed that the numerical error of the method becomes smaller even if the number of iteration steps is few, in other words, the scheme gives faster computation for nonlinear BSDEs/PDEs
than the original deep BSDE method (\cite{EHJ}) and  Fujii et al. (2019) \cite{FTT}.  

\section{New method}\label{new-method}
We propose a new method as an extension of Fujii et al. (2019) \cite{FTT} and Naito and Yamada (2020) \cite{NY}. 
The new scheme is regarded as a control variate method for solving high-dimensional nonlinear BSDEs/PDEs which is motivated by the perturbation scheme in Takahashi and Yamada (2015) \cite{TY15}. 
In the following, let us explain the proposed method. We first decompose $(Y^{\varepsilon,\alpha},Z^{\varepsilon,\alpha})$ as $Y^{\varepsilon,\alpha}=\mathcal{Y}^{1,\varepsilon}+\alpha \mathcal{Y}^{2,\varepsilon}$ and $Z^{\varepsilon,\alpha}=\mathcal{Z}^{1,\varepsilon}+\alpha \mathcal{Z}^{2,\varepsilon}$ by introducing 
\begin{align}
    -d\mathcal{Y}_t^{1,\varepsilon}&=- \mathcal{Z}_t^{1,\varepsilon} dW_t,\quad
    \mathcal{Y}_T^{1,\varepsilon}=g(X_T^\varepsilon),\label{Y1}\\
    -d\mathcal{Y}_t^{2,\varepsilon}&= f^\alpha(t,X_t^\varepsilon,\mathcal{Y}_t^{1,\varepsilon}+\alpha \mathcal{Y}_t^{2,\varepsilon},\mathcal{Z}_t^{1,\varepsilon}+ \alpha \mathcal{Z}_t^{2,\varepsilon})dt -  \mathcal{Z}_t^{2,\varepsilon} dW_t,\quad
     \mathcal{Y}_T^{2}=0. \label{Y2}
\end{align}
Here, we note that $(\mathcal{Y}^{1,\varepsilon},\mathcal{Z}^{1,\varepsilon})$ 
is the solution of a linear BSDE and that 
$(\alpha \mathcal{Y}^{2,\varepsilon},\alpha \mathcal{Z}^{2,\varepsilon})$ can be interpreted as the solution of a ``residual (nonlinear) BSDE". 

Let ${\cal U}^1$ be the solution of the linear PDE corresponding to $(\mathcal{Y}^{1,\varepsilon},\mathcal{Z}^{1,\varepsilon})$:
\begin{eqnarray}
&\partial_t {\cal U}^1(t,x)+{\cal L}^{\varepsilon}{\cal U}^1(t,x)=0, \ \ \ t<T, \\
&{\cal U}^1(T,x)=g(x).\nn
\end{eqnarray}

\subsection{Deep BSDE solver for explicitly solvable $(\mathcal{Y}^{1,\varepsilon},\mathcal{Z}^{1,\varepsilon})$}
We start with a case that $(\mathcal{Y}^{1,\varepsilon},\mathcal{Z}^{1,\varepsilon})$
is explicitly solvable as a closed-form
in order to explain our motivation of the paper. Even in this case, $(\alpha \mathcal{Y}^{2,\varepsilon},\alpha \mathcal{Z}^{2,\varepsilon})$ can not be obtained in closed-form due to the nonlinearity of the driver $f$. Hence, we apply the deep BSDE method to the residual nonlinear BSDE $(\alpha \mathcal{Y}^{2,\varepsilon},\alpha \mathcal{Z}^{2,\varepsilon})$. Then, the following will be an approximation for the target $Y_0^{\varepsilon,\alpha}$: 
\begin{eqnarray}
Y_0^{\varepsilon,\alpha} &\approx& {\cal Y}_0^{1,\varepsilon} + \alpha \widetilde{\cal Y}_0^{2,\varepsilon,(n)\ast},\label{new_decomp_solvable_case}
\end{eqnarray}
where $\widetilde{\cal Y}^{2,\varepsilon,(n)\ast}$ is obtained as a solution of the following problem based on the deep BSDE method with closed-form functions for $\mathcal{Y}^{1,\varepsilon}$ and $\mathcal{Z}^{1,\varepsilon}$: 
\begin{eqnarray}
\inf_{ \widetilde{\cal Y}^{2,\varepsilon,(n)}_0,  \widetilde{\cal Z}^{2,\varepsilon,(n)}}\Big\| \widetilde{\cal Y}_T^{2,\varepsilon,(n)} \Big\|_2^2 \label{DBSDE_optimization_prob_solvable_case}
\end{eqnarray}
subject to
\begin{align}
     \widetilde{\cal Y}_t^{2,\varepsilon,(n)}
    &= \widetilde{\cal Y}^{2,\varepsilon,(n)}_0
    -\int_0^t  f^\alpha(s,\bar{X}_s^{\varepsilon,(n)},\overline{{\cal Y}}_s^{1,\varepsilon,(n)}+\alpha \widetilde{\cal Y}_s^{2,\varepsilon,(n)},\overline{{\cal Z}}_s^{1,\varepsilon,(n)}+\alpha \widetilde{\cal Z}_s^{2,\varepsilon,(n)})ds\notag\\
    &\qquad+\int_0^t  \widetilde{\cal Z}_s^{2,\varepsilon,(n)} dW_s,
\end{align}
where 
\begin{eqnarray}
\overline{{\cal Y}}_t^{1,\varepsilon,(n)}=
{\cal U}^1(t,\bar{X}^{\varepsilon,(n)}_t), \ \ \ \ 
\overline{{\cal Z}}_t^{1,\varepsilon,(n)}=(\partial_x{\cal U}^1 \sigma^{\varepsilon})(t,\bar{X}^{\varepsilon,(n)}_t), \ \ \ \ t \in [0,T], \label{variates_solvable_case}
\end{eqnarray}
with the continuous Euler-Maruyama scheme $\bar{X}^{\varepsilon,(n)}=\{\bar{X}^{\varepsilon,(n)}_t\}_{t\geq 0}(=\bar{X}^{(n)})$
and closed-form functions
${\cal U}^1$ and
$(\partial_x{\cal U}^1 \sigma^{\varepsilon})$. 

In this case, we have the following error estimate with a small $\alpha$-effect in the residual nonlinear BSDE. 
The proof will be shown as a part of the one for Theorem \ref{main-thm-2} in the next subsection. 
Particularly, see the sentences after (\ref{error_decomposition}) and (\ref{2nd-term-error}).

\begin{thm}\label{main-thm}
There exists $C>0$ such that 
\begin{equation}
    E[|Y_0^{\varepsilon,\alpha}-\{ {\cal Y}_0^{1,\varepsilon}+\alpha \widetilde{\cal Y}_0^{2,\varepsilon,(n)} \}|^2]
    \leq \alpha^2 C \Big\{ \frac{1}{n} + \norm{ \widetilde{\mathcal{Y}}^{2,\varepsilon,(n)}_{T}}_2^2 \Big\},\label{solvable-thm}
\end{equation}
for all $\varepsilon,\alpha \in (0,1)$ and $n\geq 1$. 
\end{thm}

\subsection{General case: Deep BSDE solver for unsolvable $(\mathcal{Y}^{1,\varepsilon},\mathcal{Z}^{1,\varepsilon})$}
In most cases, $(\mathcal{Y}^{1,\varepsilon},\mathcal{Z}^{1,\varepsilon})$ is {\it unsolvable} as a closed-form, particularly it is the case when the dimension $d$ is high. In such cases, we need to approximate $(\mathcal{Y}^{1,\varepsilon},\mathcal{Z}^{1,\varepsilon})$. However, constructing tractable approximations of $\mathcal{Y}^{1,\varepsilon}_t={\cal U}^1(t,X^{\varepsilon}_t)$, $t\geq 0$, and especially $\mathcal{Z}^{1,\varepsilon}_t=(\partial_x{\cal U}^1 \sigma^{\varepsilon})(t,X^{\varepsilon}_t)$, $t \geq 0$, is not an easy task because it includes the gradient of ${\cal U}^1$. A possible solution is to use an asymptotic expansion approach with stochastic calculus. We prepare some 
notations of Malliavin calculus. Let $\mathbb{D}^\infty$ be the space of smooth Wiener functionals in the sense of Malliavin. For a {\it nondegenerate} $F \in (\mathbb{D}^\infty)^d$ and $G \in \mathbb{D}^\infty$, for a multi-index $\gamma$, there exists $H_{\gamma}(F,G) \in \mathbb{D}^\infty$ such that $E[\partial^{\gamma}\varphi(F)G]=E[\varphi(F)H_{\gamma}(F,G)]$ for all $\varphi \in C^\infty_b(\mathbb{R}^d)$. See Chapter V.8-10 in Ikeda and Watanabe (1989) \cite{IY} and Chapter 1-2 in Nualart (2006) \cite{N} for the details.  
  
First, we give approximations of ${\cal Y}^{1,\varepsilon}$ and ${\cal Z}^{1,\varepsilon}$. For $m\in \mathbb{N}$, we approximate ${\cal U}^1$ and $\partial_x {\cal U}^{1} \sigma^\varepsilon$ with asymptotic expansions 
up to the $m$-th order and Malliavin calculus, by applying or extending the methods in \cite{MTU}\cite{T2}\cite{TY12}\cite{TY15}\cite{YY19}. Let us consider $X^{t,x,\varepsilon}=\{ X^{t,x,\varepsilon}_s \}_{s\geq t}$ be the solution of 
\begin{eqnarray}
X_s^{t,x,\varepsilon}=x + \int_t^s \mu(u,X_u^{t,x,\varepsilon})du+\varepsilon  \int_t^s \sigma(u,X_u^{t,x,\varepsilon})dW_u, \ \ \ x \in \mathbb{R}^d, \  s \geq t.
\end{eqnarray}
Then the $d$-dimensional forward process $X^{t,x,\varepsilon}=(X^{t,x,\varepsilon,1},\cdots,X^{t,x,\varepsilon,d})$ can be expanded as follows: for $i=1,\cdots,d$, 
\begin{eqnarray}
X_s^{t,x,\varepsilon,i}\sim X^{t,x,0,i}_{s}+\varepsilon X^{t,x,i}_{1,s}+\varepsilon^2 X^{t,x,i}_{2,s} + \cdots \ \ \ \ \ \mbox{in} \ \ \ \mathbb{D}^\infty,
\end{eqnarray}
for some $X^{t,x,i}_{k,s} \in \mathbb{D}^\infty$, $k \in \mathbb{N}$, which are independent of $\varepsilon$ (see Watanabe (1987) \cite{W87} for example). Here, $X^{t,x,0,i}_{s}$ is the solution of $X_s^{t,x,0,i}=x + \int_t^s \mu^i(u,X_u^{t,x,0})du$, and $X^{t,x,i}_{k,s}=\frac{1}{k!}{\partial^k}/{\partial \varepsilon^k}X_s^{t,x,\varepsilon,i}|_{\varepsilon=0}$, $k \in \mathbb{N}$.

Let us define $\overline{X}^{t,x,\varepsilon}_s=X^{t,x,0}_{s}+\varepsilon X^{t,x}_{1,s}$ for $s\leq T$. The functions ${\cal U}^1$ and $\partial_x {\cal U}^{1}\sigma^\varepsilon$ are approximated by the asymptotic expansion. 
\begin{prop}\label{Prop_expansion_estimate}
Let $T>0$ and $m \in \mathbb{N}$. There is $[0,T)\times \mathbb{R}^d \times (0,1) \ni (t,x,\varepsilon) \mapsto {\cal W}^{t,x,\varepsilon,(m)}_{T} \in  \mathbb{D}^\infty$ satisfying that there exist $C(T,m)>0$ and $p(m)\geq m+1$ such that 
\begin{align}
&|{\cal U}^1(t,x)-{\cal U}^{1,(m)}(t,x) |\leq C(T,m) \varepsilon^{m+1} (T-t)^{p(m)/2}, \label{ae_error_1} 
\end{align}
for all $\varepsilon \in (0,1)$, $t<T$ and $x \in \mathbb{R}^d$, where the ${\cal U}^{1,(m)}$ is given by
\begin{align}
{\cal U}^{1,(m)}(t,x)=E[ g(\overline{X}^{t,x,\varepsilon}_T)  {\cal W}^{t,x,\varepsilon,(m)}_{T} ], \ \ \ t<T, \ x \in \mathbb{R}^d,  \label{ae_error_1-2} 
\end{align}
which satisfy ${\cal U}^{1,(m)}(t,\cdot) \in C_b^2(\mathbb{R}^d)$, $t<T$. 
Also, there is $[0,T)\times \mathbb{R}^d \times (0,1) \ni (t,x,\varepsilon) \mapsto {\cal Z}^{t,x,\varepsilon,(m)}_{T} \in  \mathbb{D}^\infty$ satisfying that there exist $K(T,m)>0$ and $q(m)\geq m$ such that 
\begin{align}
&|\partial_x {\cal U}^{1}(t,x)\sigma^\varepsilon(t,x)-{\cal V}^{1,(m)}(t,x) |\leq K(T,m) \varepsilon^{m+1} (T-t)^{q(m)/2}, \label{ae_error_2} 
\end{align}
for all $\varepsilon \in (0,1)$, $t<T$ and $x \in \mathbb{R}^d$, where the ${\cal V}^{1,(m)}$ is given by
\begin{align}
{\cal V}^{1,(m)}(t,x)=E[ g(\overline{X}^{t,x,\varepsilon}_T)  {\cal Z}^{t,x,\varepsilon,(m)}_{T} ], \ \ \ t<T, \ x \in \mathbb{R}^d, \label{ae_error_2-2} 
\end{align} 
which satisfy ${\cal V}^{1,(m)}(t,\cdot) \in C_b^1(\mathbb{R}^d)$, $t<T$. 
\end{prop}
\noindent
{\it Proof}. See Appendix \ref{Proof_Prop_expansion_estimate}. \\

For example, the stochastic weight ${\cal W}^{t,x,\varepsilon,(m)}_{T}$ has the representation in general: 
\begin{align}
& {\cal W}^{t,x,\varepsilon,(m)}_{T} \nn\\
=&1+\sum_{j=1}^m \varepsilon^{j} \sum_{k=1}^j \sum_{\beta_1+\cdots+\beta_k=j,\beta_i\geq 1}\sum_{\gamma^{(k)}=(\gamma_1,\cdots,\gamma_k)\in\{1,\cdots,d \}^k}\frac{1}{k!} H_{\gamma^{(k)}}(X^{t,x}_{1,T},\prod_{\ell=1}^k X^{t,x,\gamma_\ell}_{{\beta_\ell}+1,T} ).
\end{align}
See Section 2.2 in Takahashi and Yamada (2012) \cite{TY12} and Section 6.1 in Takahashi (2015) \cite{T2} for more details. 
The functions ${\cal U}^{1,(m)}$ and ${\cal V}^{1,(m)}$ have more explicit representation. Actually when $m=1$, ${\cal U}^{1,(1)}$ and ${\cal V}^{1,(1)}$ have the following forms which are easily computed by taking advantage of the fact that $\overline{X}_T^{t,x,\varepsilon}$ (and $X_{1,T}^{t,x}$) is a Gaussian random variable. In particular, the representation ${\cal V}^{1,(1)}$, the multidimensional expansion of $\partial_x {\cal U}^{1}\sigma^\varepsilon$ is new, which is an extension of \cite{MTU}\cite{YY19}. 

\begin{prop}\label{first_order_ae_U}
For $t<T$, $x \in \mathbb{R}^d$, 
\begin{align}
  &{\cal U}^{1,(1)}(t,x)
  =E[g(\overline{X}_T^{t,x,\varepsilon})] \label{cal-U}\\
  &+ \varepsilon
  \sum_{i_1,i_2,i_3,j_1=1}^d \sum_{k_1,k_2=1}^d
     E[ g(\overline{X}_T^{t,x,\varepsilon}) H_{(i_1,i_2,i_3)}(X_{1,T}^{t,x},1) ] \ C_{i_1,i_2,i_3,j_1}^{(1),k_1,k_2}(t,T,x)\nn  \\
  &+\varepsilon \sum_{i_1,i_2,i_3,j_1,j_2=1}^d \sum_{k_1,k_2=1}^d E[ g(\overline{X}_T^{t,x,\varepsilon})  H_{(i_1,i_2,i_3)}(X_{1,T}^{t,x},1) ] \ C_{i_1,i_2,i_3,j_1,j_2}^{(2),k_1,k_2}(t,T,x) \nn \\
  &+\varepsilon \frac{1}{2} \sum_{i_1,j_1,j_2=1}^d \sum_{k_1,k_2=1}^d  E[ g(\overline{X}_T^{t,x,\varepsilon})  H_{(i_1)}(X_{1,T}^{t,x},1) ] \mathrm{1}_{k_1=k_2}C_{i_1,j_1,j_2}^{(3),k_1,k_2}(t,T,x), \nn
\end{align}
{\footnotesize{\begin{align}
  &{\cal V}^{1,(1)}(t,x)= \sum_{i_1=1}^d E[g(\overline{X}_T^{t,x,\varepsilon}) H_{(i_1)}(X_{1,T}^{t,x},1)] [J_{t\to T}^{0,x}]^{i_1}\sigma(t,x)
  \label{cal V-1} \\
  &+ \varepsilon \sum_{i_1,i_2,i_3,i_4,j_1=1}^d \sum_{k_1,k_2=1}^d
E[ g(\overline{X}_T^{t,x,\varepsilon}) H_{(i_1,i_2,i_3,i_4)}(X_{1,T}^{t,x},1)] \ [J_{t\to T}^{0,x}]^{i_1} C_{i_2,i_3,i_4,j_1}^{(1),k_1,k_2}(t,T,x) \sigma(t,x)\nn \\
  &+\varepsilon \sum_{i_1,i_2,i_3,i_4,j_1,j_2=1}^d \sum_{k_1,k_2=1}^d E[g(\overline{X}_T^{t,x,\varepsilon}) H_{(i_1,i_2,i_3,i_4)}(X_{1,T}^{t,x},1)] \ [J_{t\to T}^{0,x}]^{i_1} C_{i_2,i_3,i_4,j_1,j_2}^{(2),k_1,k_2}(t,T,x)\sigma(t,x)\nn \\
  &+\varepsilon \frac{1}{2} \sum_{i_1,j_1,j_2=1}^d \sum_{k_1,k_2=1}^d E[ g(\overline{X}_T^{t,x,\varepsilon}) H_{(i_1,i_2)}(X_{1,T}^{t,x},1)] [J_{t\to T}^{0,x}]^{i_1} \mathrm{1}_{k_1=k_2}C_{i_2,j_1,j_2}^{(3),k_1,k_2}(t,T,x)\sigma(t,x)\nn\\
  &+ \varepsilon \sum_{i_1,i_2,j_1,j_2=1}^d \sum_{k_1=1}^d E[g(\overline{X}_T^{t,x,\varepsilon}) H_{(i_1,i_2)}(X_{1,T}^{t,x},1)] [J_{t\to T}^{0,x}]_{j_1}^{i_1} C_{i_2,j_1,j_2}^{(4),k_1}(t,T,x) \sigma(t,x)\nn \\
  &+ \varepsilon \sum_{i_1,i_2,j_1=1}^d \sum_{k_1=1}^d E[g(\overline{X}_T^{t,x,\varepsilon}) H_{(i_1,i_2)}(X_{1,T}^{t,x},1)] [J_{t\to T}^{0,x}]_{j_1}^{i_1} C_{i_2,j_1}^{(5),k_1}(t,T,x) \sigma(t,x),\nn
\end{align}}}
and 
{\footnotesize{\begin{align*}
C_{i_1,i_2,i_3,j_1}^{(1),k_1,k_2}(t,T,x)=&\int_t^T \int_t^{t_1} a^{i_3}_{k_2}(t,t_2,t_1,x) a^{i_2}_{k_1}(t,t_1,T,x) b^{i_1,j_1}_{k_1}(t,t_1,T,x) a^{j_1}_{k_2}(t,t_2,t_1,x)  dt_2dt_1,\\
C_{i_1,i_2,i_3,j_1,j_2}^{(2),k_1,k_2}(t,T,x)
=&\int_t^T \int_t^{t_1} \int_t^{t_2}
  a^{i_3}_{k_1}(t,t_3,t_2,x)
  a^{i_2}_{k_2}(t,t_2,t_1,x)\\
  & \hspace{4em} c^{i_1,j_1,j_2}(t,t_1,T,x)
  a^{j_1}_{k_2}(t,t_2,t_1,x)
  a^{j_2}_{k_1}(t,t_3,t_1,x)
  dt_3
  dt_2
  dt_1,\\
C_{i_1,j_1,j_2}^{(3),k_1,k_2}(t,T,x)=&\int_t^T \int_t^{t_1} c^{i_1,j_1,j_2}(t,t_1,T,x) a^{j_2}_{k_2}(t,t_2,t_1,x) a^{j_1}_{k_1}(t,t_2,t_1,x) dt_2 dt_1,
\end{align*}}}
{\footnotesize{\begin{align*}
  C_{i_1,j_1,j_2}^{(4),k_1}(t,T,x)=&\int_t^T \int_t^{t_1} a^{i_1}_{k_1}(t,t_2,T,x) [\partial^2 \mu(t_1,X_{t_1}^{t,x,0})]^{j_1}_{j_2} a^{j_2}_{k_1}(t,t_2,t_1,x) dt_2 dt_1,\\
  C_{i_1,j_1}^{(5),k_1}(t,T,x)=&\int_t^T a_{k_1}^{i_1}(t,t_1,T,x)\partial_{j_1} \sigma_{k_1}(t_1,X_{t_1}^{t,x,0})dt_1,
\end{align*}}}
with
{\footnotesize{\begin{align*}
a^i_k(t,s,u,x)&:=\sum_{j_1,j_2=1}^d [J_{t\to u}^{0,x}]^i_{j_1} [(J_{t \to s}^{0,x})^{-1}]^{j_1}_{j_2} \sigma_k^{j_2}(s,X_s^{t,x,0}),\\
b^{i,j_3}_k(t,s,u,x)&:=\sum_{j_1,j_2=1}^d [J_{t\to u}^{0,x}]^i_{j_1} [(J_{t \to s}^{0,x})^{-1}]^{j_1}_{j_2} \partial_{j_3} \sigma_k^{j_2}(s,X_s^{t,x,0}),\\
c^{i,j_3,j_4}(t,s,u,x)&:=\sum_{j_1,j_2=1}^d [J_{t\to u}^{0,x}]^i_{j_1} [(J_{t \to s}^{0,x})^{-1}]^{j_1}_{j_2} [\partial^2 \mu^{j_2}(s,X_s^{t,x,0})]^{j_3}_{j_4}.
\end{align*}}}
Here, $[ \ \cdot \ ]^i_j$ is an entry in $i$-th row and $j$-th column of a matrix, $\partial_{j} \varphi(\cdot)$ is an $j$-th element of $\partial \varphi(\cdot) = \left[\pdv*{\varphi}{x_i} (\cdot)\right]_{1\leq i\leq d}$ and $[\partial^2\varphi(\cdot)]^i_j=\pdv[2]{\varphi(\cdot)}{x_i}{x_j}$ 
($1\leq i,j \leq d$)
is used. 
\end{prop}
\noindent
{\it Proof}. See Appendix \ref{Proof_Prop_expansion_order1}. \\ 
 
Using ${\cal U}^{1,(m)}$ and ${\cal V}^{1,(m)}$, we define 
 \begin{eqnarray}
\overline{{\cal Y}}_t^{1,\varepsilon,(m)}={\cal U}^{1,(m)}(t,X_t^\varepsilon), \ \ \ \overline{{\cal Z}}_t^{1,\varepsilon,(m)}={\cal V}^{1,(m)}(t,X_t^\varepsilon), \ \ \ t\geq 0. 
\end{eqnarray}

Furthermore, we compute $ {\cal Y}^{2,\varepsilon}$ and $ {\cal Z}^{2,\varepsilon}$ numerically by the deep BSDE method by solving
\begin{eqnarray}
\inf_{ {\cal Y}^{2,\varepsilon,(m,n)}_0,  {\cal Z}^{2,\varepsilon,(m,n)}}\Big\| {\cal Y}_T^{2,\varepsilon,(m.n)} \Big\|_2^2 \label{DBSDE_optimization_prob}
\end{eqnarray}
subject to
\begin{align}
     {\cal Y}_t^{2,\varepsilon,(m,n)}
    &= {\cal Y}^{2,\varepsilon,(m,n)}_0
    -\int_0^t  f(s,\bar{X}_s^{\varepsilon,(n)},\overline{{\cal Y}}_s^{1,\varepsilon,(m,n)}+\alpha {\cal Y}_s^{2,\varepsilon,(m,n)},\overline{{\cal Z}}_s^{1,\varepsilon,(m,n)}+\alpha {\cal Z}_s^{2,\varepsilon,(m,n)})ds\notag\\
    &\qquad+\int_0^t  {\cal Z}_s^{2,\varepsilon,(m,n)} dW_s,
\end{align}
where 
\begin{eqnarray}
\overline{{\cal Y}}_t^{1,\varepsilon,(m,n)}={\cal U}^{1,(m)}(t,\bar{X}^{\varepsilon,(n)}_t), \ \ \ \ 
\overline{{\cal Z}}_t^{1,\varepsilon,(m,n)}={\cal V}^{1,(m)}(t,\bar{X}^{\varepsilon,(n)}_t), \ \ \ \ t \in [0,T], \label{linear-ae}
\end{eqnarray}
with the continuous Euler-Maruyama scheme $\bar{X}^{\varepsilon,(n)}=\{\bar{X}^{\varepsilon,(n)}_t\}_{t\geq 0}(=\bar{X}^{(n)})$. \\

We have the main theoretical result in this paper as follows. 
\begin{thm}\label{main-thm-2}
There exists $C>0$ such that
\begin{equation}
    E[|Y_0^{\varepsilon,\alpha}-\{ \overline{{\cal Y}}_0^{1,\varepsilon,(m)}+\alpha {\cal Y}_0^{2,\varepsilon,(m,n)} \}|^2]
    \leq C \varepsilon^{2(m+1)}
    +\alpha^2 C \Big\{\varepsilon^{2(m+1)} +\frac{1}{n} + \norm{\mathcal{Y}^{2,\varepsilon,(m,n)}_{T}}_2^2 \Big\},\label{unsolvable-main-thm}
\end{equation}
for all $\varepsilon,\alpha \in (0,1)$ and $n\geq 1$. 
\end{thm}

\noindent
{\it Proof}. \ 
  In the proof, we use a generic constant $C>0$ which varies from line to line.
  Let $( {\cal Y}^{2,\varepsilon,(m)},  {\cal Z}^{2,\varepsilon,(m)})$ be the solution of the following BSDE:
  \begin{align}
     {\cal Y}_t^{2,\varepsilon,(m)}
    &=\int_t^T  f(s,X_s^\varepsilon,\overline{{\cal Y}}_s^{1,\varepsilon,(m)}+\alpha {\cal Y}_s^{2,\varepsilon,(m)},\overline{{\cal Z}}_s^{1,\varepsilon,(m)}+\alpha {\cal Z}_s^{2,\varepsilon,(m)})ds
    -\int_t^T  {\cal Z}_s^{2,\varepsilon,(m)} dW_s.
  \end{align}
  Then we have 
  \begin{align}
      &\quad E[|Y_0^{\varepsilon,\alpha}-\{ \overline{{\cal Y}}_0^{1,\varepsilon,(m)}+\alpha {\cal Y}_0^{2,\varepsilon,(m,n)} \}|^2]\nn\\
      &=E[|Y_0^{\varepsilon,\alpha}-\{ \overline{{\cal Y}}_0^{1,\varepsilon,(m)}+\alpha {\cal Y}_0^{2,\varepsilon,(m)} \} + \alpha {\cal Y}_0^{2,\varepsilon,(m)} - \alpha {\cal Y}_0^{2,\varepsilon,(m,n)}|^2]\nn\\
      &\leq 
      CE[|Y_0^{\varepsilon,\alpha}-\{ \overline{{\cal Y}}_0^{1,\varepsilon,(m)}+\alpha {\cal Y}_0^{2,\varepsilon,(m)} \}|^2]
      +\alpha^2 CE[| {\cal Y}_0^{2,\varepsilon,(m)} - {\cal Y}_0^{2,\varepsilon,(m,n)}|^2]. \label{error_decomposition}
  \end{align}
First, we estimate the term   
$E[|Y_0^{\varepsilon,\alpha}-\{ \overline{{\cal Y}}_0^{1,\varepsilon,(m)}+\alpha {\cal Y}_0^{2,\varepsilon,(m)} \}|^2]$. We note that this term becomes null in 
(\ref{solvable-thm}), i.e.
the error estimate of Theorem \ref{main-thm} for
the case that $(\mathcal{Y}^{1,\varepsilon},\mathcal{Z}^{1,\varepsilon})$
is explicitly solvable as a closed-form.

  Since we have
  \begin{eqnarray}
  Y_0^{\varepsilon,\alpha}=E_{X_0}[g(X_T^\varepsilon)]+\alpha E_{X_0}[\int_0^T  f(s,X_s^\varepsilon,Y_s^{\varepsilon,\alpha},Z_s^{\varepsilon,\alpha})ds]
  \end{eqnarray}
  and
  \begin{align}
    \overline{{\cal Y}}_0^{1,\varepsilon,(m)}
    &=E_{X_0}[g(\overline{X}^{0,\cdot,\varepsilon}_T)\{ 1 + {\cal W}^{0,\cdot,\varepsilon,(m)}_{T} \}],\\
     {\cal Y}^{2,\varepsilon,(m)}_0
    &=E_{X_0}[\int_0^T  f(s,X_s^\varepsilon,\overline{{\cal Y}}_s^{1,\varepsilon,(m)}+\alpha {\cal Y}_s^{2,\varepsilon,(m)},\overline{{\cal Z}}_s^{1,\varepsilon,(m)}+\alpha {\cal Z}_s^{2,\varepsilon,(m)})ds],
  \end{align}
  it holds that
  \begin{align}
      &\quad E[|Y_0^{\varepsilon,\alpha}-\{ \overline{{\cal Y}}_0^{1,\varepsilon,(m)}+\alpha {\cal Y}_0^{2,\varepsilon,(m)} \}|^2] \notag\\
      &\leq C E[|E_{X_0}[g(X_T^\varepsilon)]-E_{X_0}[g(\overline{X}^{0,\cdot,\varepsilon}_T) {\cal W}^{0,\cdot,\varepsilon,(m)}_{T} ]|^2] \notag\\
      &\quad+C E\Big[ \Big|E_{X_0}[\int_0^T \alpha f(s,X_s^\varepsilon,Y_s^{\varepsilon,\alpha},Z_s^{\varepsilon,\alpha})ds]\nn\\
      &\qquad\qquad-E_{X_0}[\int_0^T \alpha f(s,X_s^\varepsilon,\overline{{\cal Y}}_s^{1,\varepsilon,(m)}+\alpha {\cal Y}_s^{2,\varepsilon,(m)},\overline{{\cal Z}}_s^{1,\varepsilon,(m)}+\alpha {\cal Z}_s^{2,\varepsilon,(m)})ds] \Big|^2 \Big]\notag\\
      &\leq C \varepsilon^{2(m+1)} \nn\\
      &+ C \alpha^2 C_{\mathrm{Lip}}[f]^2 \int_0^TE[ | {\cal Y}_s^{1,\varepsilon}-\overline{{\cal Y}}_s^{1,\varepsilon,(m)} |^2 ] ds
      + C \alpha^2 C_{\mathrm{Lip}}[f]^2 \int_0^TE[ | {\cal Z}_s^{1,\varepsilon}-\overline{{\cal Z}}_s^{1,\varepsilon,(m)} |^2 ] ds\notag\\
      &+ C \alpha^2 C_{\mathrm{Lip}}[f]^2 \int_0^TE[ | \alpha {\cal Y}_s^{2,\varepsilon} - \alpha {\cal Y}_s^{2,\varepsilon,(m)} |^2 ] ds
      + C \alpha^2 C_{\mathrm{Lip}}[f]^2 \int_0^TE[ | \alpha {\cal Z}_s^{2,\varepsilon} - \alpha {\cal Z}_s^{2,\varepsilon,(m)} |^2 ] ds. 
  \end{align}
  Here, the estimates
  \begin{align}
  &\int_0^TE[ | {\cal Y}_s^{1,\varepsilon}-\overline{{\cal Y}}_s^{1,\varepsilon,(m)} |^2 ] ds \leq C \varepsilon^{2(m+1)}, \label{est_ae-0}\\
  &\int_0^TE[ | {\cal Z}_s^{1,\varepsilon}-\overline{{\cal Z}}_s^{1,\varepsilon,(m)} |^2 ] ds\leq C \varepsilon^{2(m+1)}, \label{est_ae}
  \end{align}
  are obtained by (\ref{ae_error_1}) and (\ref{ae_error_2}). Also, by Theorem~{4.2.3} in Zhang (2017) \cite{Z}, we have 
  \begin{align}
    & \int_0^TE[ | \alpha {\cal Y}_s^{2,\varepsilon} - \alpha {\cal Y}_s^{2,\varepsilon,(m)} |^2 ] ds
      +  \int_0^TE[ | \alpha {\cal Z}_s^{2,\varepsilon} - \alpha {\cal Z}_s^{2,\varepsilon,(m)} |^2 ] ds \nn\\
    &\leq C E[
      \int_0^T |\alpha f(s,X_s^\varepsilon,\mathcal{Y}_s^{1,\varepsilon}+\alpha \mathcal{Y}_s^{2,\varepsilon},\mathcal{Z}_s^{1,\varepsilon}+ \alpha \mathcal{Z}_s^{2,\varepsilon}) \nn\\
      &\qquad\qquad\qquad- \alpha f(s,X_s^\varepsilon,\overline{{\cal Y}}_s^{1,\varepsilon,(m)}+\alpha \mathcal{Y}_s^{2,\varepsilon},\overline{{\cal Z}}_s^{1,\varepsilon,(m)}+\alpha \mathcal{Z}_s^{2,\varepsilon})|^2 ds
    ]\nn\\
    &\leq C \alpha^2 C_{\mathrm{Lip}}[f]^2
    \{
      \int_0^TE[ | {\cal Y}_s^{1,\varepsilon}-\overline{{\cal Y}}_s^{1,\varepsilon,(m)} |^2 ] ds
      + \int_0^TE[ | {\cal Z}_s^{1,\varepsilon}-\overline{{\cal Z}}_s^{1,\varepsilon,(m)} |^2 ] ds
    \}\nn\\
    &\leq C \alpha^2 \varepsilon^{2(m+1)},
  \end{align}
  where the estimates (\ref{est_ae-0}) and (\ref{est_ae}) are applied in the last inequality. 
  Therefore, we get 
  \begin{eqnarray}
  E[|Y_0^{\varepsilon,\alpha}-\{ \overline{{\cal Y}}_0^{1,\varepsilon,(m)}+\alpha {\cal Y}_0^{2,\varepsilon,(m)} \}|^2] \leq C \varepsilon^{2(m+1)}+C \alpha^2 \varepsilon^{2(m+1)}. 
  \end{eqnarray}

  Next, we estimate
  \begin{eqnarray}
  E[| {\cal Y}_0^{2,\varepsilon,(m)}- {\cal Y}_0^{2,\varepsilon,(m,n)}|^2]
  \label{2nd-term-error}
  \end{eqnarray}
  in (\ref{error_decomposition}). 
  We note that only this term appears in (\ref{solvable-thm}), i.e.
  the error estimate of Theorem \ref{main-thm} 
  for the case that $(\mathcal{Y}^{1,\varepsilon},\mathcal{Z}^{1,\varepsilon})$
  is explicitly solvable as a closed-form.
  
  Since we have
  \begin{align}
    & {\cal Y}_0^{2,\varepsilon,(m)}- {\cal Y}^{2,\varepsilon,(m,n)}_0\nn\\
    &=\int_0^T  f(s,X_s^\varepsilon,\overline{{\cal Y}}_s^{1,\varepsilon,(m)}+\alpha {\cal Y}_s^{2,\varepsilon,(m)},\overline{{\cal Z}}_s^{1,\varepsilon,(m)}+\alpha {\cal Z}_s^{2,\varepsilon,(m)})ds-\int_0^T  {\cal Z}_s^{2,\varepsilon,(m)} dW_s \nonumber\\
    &\quad - {\cal Y}_T^{2,\varepsilon,(m,n)}-\int_0^T  f(s,\bar{X}_s^{\varepsilon,(n)},\overline{{\cal Y}}_s^{1,\varepsilon,(m,n)}+\alpha {\cal Y}_s^{2,\varepsilon,(m,n)},\overline{{\cal Z}}_s^{1,\varepsilon,(m,n)}+\alpha {\cal Z}_s^{2,\varepsilon,(m,n)})ds \nn\\
    &\quad+\int_0^T  {\cal Z}_s^{2,\varepsilon,(m,n)} dW_s,
  \end{align}
 the upper bound of $E[| {\cal Y}_0^{2,\varepsilon,(m)}- {\cal Y}_0^{2,\varepsilon,(m,n)}|^2] $ can be decomposed as 
  \begin{align}
      &E[| {\cal Y}_0^{2,\varepsilon,(m)}- {\cal Y}_0^{2,\varepsilon,(m,n)}|^2] \notag\\
      &\leq \| \mathcal{Y}^{2,\varepsilon,(m,n)}_{T}\|_2^2
      +\int_0^T E[| {\cal Z}_s^{2,\varepsilon,(m)} -  {\cal Z}_s^{2,\varepsilon,(m,n)}|^2] ds \nn\\
      &\quad+C_{\text{Lip}}[f]^2 \times \Big\{
      E[
        \int_0^T |X_s^\varepsilon - \bar{X}_s^{\varepsilon,(n)}|^2 ds\nn
        \\
        &\quad+\int_0^T |\overline{{\cal Y}}_s^{1,\varepsilon,(m)}- \overline{{\cal Y}}_s^{1,\varepsilon,(m,n)}|^2 ds
        +\int_0^T |\alpha {\cal Y}_s^{2,\varepsilon,(m)} - \alpha {\cal Y}_s^{2,\varepsilon,(m,n)}|^2 ds\nn\\
        &\quad+\int_0^T |\overline{{\cal Z}}_s^{1,\varepsilon,(m)} - \overline{{\cal Z}}_s^{1,\varepsilon,(m,n)}|^2 ds
        +\int_0^T |\alpha {\cal Z}_s^{2,\varepsilon,(m)} - \alpha {\cal Z}_s^{2,\varepsilon,(m,n)}|^2 ds
      ]
      \Big\}.
  \end{align}
Then, the following holds:
\begin{align}
&\int_0^T E[|\overline{{\cal Y}}_s^{1,\varepsilon,(m)}- \overline{{\cal Y}}_s^{1,\varepsilon,(m,n)}|^2] ds \leq C E[
\int_0^T |X_s^\varepsilon - \bar{X}_s^{\varepsilon,(n)}|^2 ds], \\
&\int_0^T E[|\overline{{\cal Z}}_s^{1,\varepsilon,(m)}- \overline{{\cal Z}}_s^{1,\varepsilon,(m,n)}|^2] ds \leq C E[
\int_0^T |X_s^\varepsilon - \bar{X}_s^{\varepsilon,(n)}|^2 ds],            
\end{align} 
 since for all $t<T$, ${\cal U}^{1,(m)}(t,\cdot)$ and ${\cal V}^{1,(m)}(t,\cdot)$ are in $C^2_b$ and $C^1_b$, respectively. Thus, we have
   \begin{align}
      &  E[| {\cal Y}_0^{2,\varepsilon,(m)}- {\cal Y}_0^{2,\varepsilon,(m,n)}|^2] \leq \|\mathcal{Y}^{2,\varepsilon,(m,n)}_{T}\|_2^2+C \times
      \Big\{
        \sup_{t\in[0,T]} (E[|X_t-\bar{X}^{\varepsilon,(n)}_t|^2] \nn\\
      &\quad \ \ \ \ \ \ \ \ \ \ +E[|{\cal Y}_t^{2,\varepsilon,(m)}-{\cal Y}_t^{2,\varepsilon,(m,n)}|^2])
        +\int_0^T E[|{\cal Z}_s^{2,\varepsilon,(m)}-{\cal Z}_s^{2,\varepsilon,(m,n)}|^2]ds
      \Big\}.
  \end{align}
  By Theorem 1 of Han and Long (2020) \cite{HL}, it holds that 
  \begin{align}
      &\sup_{t\in[0,T]} (E[|X_t^\varepsilon-\bar{X}^{\varepsilon,(n)}_t|^2]+E[|{\cal Y}_t^{2,\varepsilon,(m)}-{\cal Y}_t^{2,\varepsilon,(m,n)}|^2]) + \int_0^T E[|{\cal Z}_s^{2,\varepsilon,(m)}-{\cal Z}_s^{2,\varepsilon,(m,n)}|^2]ds \nn\\
      &\leq C \Big\{\frac{T}{n} + \norm{\mathcal{Y}^{2,\varepsilon,(m,n)}_{T}}_2^2\Big\}. 
  \end{align}
  Therefore, we get
  \begin{equation}
    E[| {\cal Y}_0^{2,\varepsilon,(m)}- {\cal Y}_0^{2,\varepsilon,(m,n)}|^2]\\
    \leq \frac{C}{n} + C \norm{\mathcal{Y}^{2,\varepsilon,(m,n)}_{T}}_2^2,
  \end{equation}
  and the assertion is obtained as: 
   \begin{align}
      &\quad E[|Y_0^{\varepsilon,\alpha}-\{ \overline{{\cal Y}}_0^{1,\varepsilon,(m)}+\alpha {\cal Y}_0^{2,\varepsilon,(m,n)} \}|^2]\\
      &\leq CE[|Y_0^{\varepsilon,\alpha}-\{ \overline{{\cal Y}}_0^{1,\varepsilon,(m)}+\alpha {\cal Y}_0^{2,\varepsilon,(m)} \}|^2]+\alpha^2 CE[| {\cal Y}_0^{2,\varepsilon,(m)} - {\cal Y}_0^{2,\varepsilon,(m,n)}|^2]\nn\\
      &\leq C \varepsilon^{2(m+1)} + C \alpha^2\varepsilon^{2(m+1)} + \alpha^2 C \Big\{ \frac{1}{n} + \norm{\mathcal{Y}^{2,\varepsilon,(m,n)}_{T}}_2^2 \Big\}. \ \ \ \ \ \ \Box \nn
  \end{align} 

By the theorem above, it holds that  
\begin{eqnarray}
Y_0^{\varepsilon,\alpha} &\approx& \overline{{\cal Y}}_0^{1,\varepsilon,(m)} + \alpha {\cal Y}_0^{2,\varepsilon,(m,n)\ast},\label{new_decomp}
\end{eqnarray}
where ${\cal Y}_0^{2,\varepsilon,(m,n)\ast}$ is obtained by solving (\ref{DBSDE_optimization_prob}) with Deep BSDE method. The process $\overline{{\cal Y}}^{1,\varepsilon,(m,n)}$ and $\overline{{\cal Z}}^{1,\varepsilon,(m,n)}$ work as control variates for the nonlinear BSDE. 

Here, let us briefly make comments on comparison of the theoretical error estimates of our proposed method, namely (\ref{solvable-thm}) in Theorem \ref{main-thm} 
for the explicitly solvable $(\mathcal{Y}^{1,\varepsilon},\mathcal{Z}^{1,\varepsilon})$ case
and 
(\ref{unsolvable-main-thm}) in Theorem \ref{main-thm-2}
for the unsolvable $(\mathcal{Y}^{1,\varepsilon},\mathcal{Z}^{1,\varepsilon})$ case
with the one provided by Han and Long (2020)
for the method of Weinan E et al (2017), i.e. 
(\ref{han-long-thm}) in Theorem \ref{HL-thm1}.
Given the number of discretized time steps $n$ for Euler-Maruyama scheme,
those are relisted below:
\begin{itemize}
\item \mbox{Proposed method (for the solvable $(\mathcal{Y}^{1,\varepsilon},\mathcal{Z}^{1,\varepsilon})$ case):}
\begin{eqnarray*}
E[|Y_0^{\varepsilon,\alpha}-\{ {\cal Y}_0^{1,\varepsilon}+\alpha \widetilde{\cal Y}_0^{2,\varepsilon,(n)} \}|^2]
	\leq C \alpha^2  \frac{1}{n} 
	+  C \alpha^2 \norm{ \widetilde{\mathcal{Y}}^{2,\varepsilon,(n)}_{T}}_2^2
	\ \ \ 
	(\ref{solvable-thm})
\end{eqnarray*}
\item \mbox{Proposed method (for the unsolvable $(\mathcal{Y}^{1,\varepsilon},\mathcal{Z}^{1,\varepsilon})$ case):}
\begin{eqnarray*}
&&E[|Y_0^{\varepsilon,\alpha}-\{ \overline{{\cal Y}}_0^{1,\varepsilon,(m)}+\alpha {\cal Y}_0^{2,\varepsilon,(m,n)} \}|^2]\\
&& \ \ \ \ \ \ \ \leq 
	(C \varepsilon^{2(m+1)}+C \alpha^2 \varepsilon^{2(m+1)}) 
	+C \alpha^2 \frac{1}{n} 
	+ C \alpha^2 \norm{\mathcal{Y}^{2,\varepsilon,(m,n)}_{T}}_2^2
	\ \ (\ref{unsolvable-main-thm})
\end{eqnarray*}	
\item \mbox{Method of Weinan E et al. (2017) (error estimate by Han and Long (2020)):}
\begin{eqnarray*}	
 E[|Y_0^{\varepsilon,\alpha} - Y_0^{\varepsilon,\alpha,(n)}|^2]
\leq C \frac{1}{n} + C \norm{g(\bar{X}_T^{\varepsilon,(n)}) - {Y}_T^{\varepsilon,\alpha,(n)}}_2^2.\ \ (\ref{han-long-thm})
\end{eqnarray*}
\end{itemize}
Thanks to the following advantages of our proposed method,
we can see that it works better as a new Deep BSDE solver, more precisely,
its errors are expected to be smaller:
\begin{itemize}
\item
(i)
Decomposition into a ``dominant" linear PDE with original terminal $g$
and a ``small" nonlinear PDE with zero terminal, i.e.
\begin{align*}
	u(0,x)=\underset{\mbox{\tiny{``dominant" linear PDE part}}}{{\cal U}^1(0,x)}+\underset{\mbox{\tiny{``small" nonlinear PDE part}}}{{\cal U}^2(0,x)}, \ \ \ \ x\in \mathbb{R}^d. \label{decomposition}
\end{align*}
and an application of Deep BSDE solver only to the 
``small" nonlinear PDE.

(ii) Closed form solutions/approximations 
for the linear PDE,
which also work as control variates for the driver of the nonlinear PDE.

Thanks to (i) and (ii), we can obtain the term
$C \alpha^2 \norm{ \widetilde{\mathcal{Y}}^{2,\varepsilon,(n)}_{T}}_2^2$
in the error bound,
rather than 
$C \norm{g(\bar{X}_T^{\varepsilon,(n)}) - {Y}_T^{\varepsilon,\alpha,(n)}}_2^2$.

Moreover, we note that our method enjoys the effects of
a small parameter in the nonlinear driver $\alpha\in (0,1)$ for this term,
as well as for the discretization error term caused by Euler-Maruyama scheme, 
which is given as $C \alpha^2  \frac{1}{n}$
rather than $C \frac{1}{n}$.

\item
Regarding
the unsolvable $(\mathcal{Y}^{1,\varepsilon},\mathcal{Z}^{1,\varepsilon})$ case,
our asymptotic expansions with respect to a small parameter $\varepsilon\in (0,1)$
in the diffusion coefficient
enable us to obtain closed form approximations 
$\overline{{\cal Y}}_t^{1,\varepsilon,(m)}={\cal U}^{1,(m)}(t,X_t^\varepsilon)$
and 
$\overline{{\cal Z}}_t^{1,\varepsilon,(m)}={\cal V}^{1,(m)}(t,X_t^\varepsilon)$: 
Particularly, 
in (\ref{unsolvable-main-thm})
the coefficients
$C \varepsilon^{2(m+1)}$
and
$C \alpha^2 \varepsilon^{2(m+1)}$
are associated with errors of the approximations for terminal $g$
and driver $\alpha f$, respectively.
\end{itemize}

We will check the effectiveness of the new method by numerical experiments in the next section.   

\begin{rem}\label{remark_1}
We give an important remark on the new method. While the proposed scheme provides a fine result, we can further improve it by replacing our approximation for the linear part $\overline{{\cal Y}}_0^{1,\varepsilon,(m)}$ in the decomposition (\ref{new_decomp}) with 
the methods of \cite{TYo}\cite{TY16}\cite{Y19}\cite{NYmcma}\\ \cite{OY}\cite{IY}. 

For example, based on Takahashi and Yamada (2016) \cite{TY16}, the following result will be an improvement of the proposed scheme. Let $t_i=T(1-(1-i/n_0)^{\gamma})$, $i=0,1,\cdots,n$, with a parameter $\gamma >0$, and $\bar{X}_{t_i}^{0,x,\varepsilon,(n)}=\overline{X}_{t_i}^{t_{i-1},\bar{X}_{t_{i-1}}^{0,x,\varepsilon,(n)},\varepsilon}$, $i=1,\cdots,n$. Define 
\begin{eqnarray}
\widehat{{\cal Y}}_0^{1,\varepsilon,(m,n_0)}=E[ g(\bar{X}_T^{0,x,\varepsilon,(n_0)}) \prod_{i=1}^{n_0} {\cal W}_{t_i}^{t_{i-1},\bar{X}_{t_{i-1}}^{0,x,\varepsilon,(n_0)},\varepsilon}  ]|_{x=X_0},
\end{eqnarray}
and consider the quantity 
\begin{eqnarray}
\widehat{{\cal Y}}_0^{1,\varepsilon,(m,n_0)} + \alpha {\cal Y}_0^{2,\varepsilon,(m,n)\ast}, \label{improved_approx}
\end{eqnarray}
where ${\cal Y}_0^{2,\varepsilon,(m,n)\ast}$ is the same as in (\ref{new_decomp}). Then, (\ref{improved_approx}) will be the improved approximation, as   
\begin{eqnarray}
Y_0^{\varepsilon,\alpha} &\approx& \widehat{{\cal Y}}_0^{1,\varepsilon,(m,n_0)} + \alpha {\cal Y}_0^{2,\varepsilon,(m,n)\ast}, 
\end{eqnarray}
in the following sense.
\begin{cor}\label{improved_method1}
There exist $C>0$ and $r(m)>0$ such that
\begin{equation}
    E[|Y_0^{\varepsilon,\alpha}-\{ \widehat{{\cal Y}}_0^{1,\varepsilon,(m,n_0)}+\alpha {\cal Y}_0^{2,\varepsilon,(m,n)} \}|^2]
    \leq C \frac{\varepsilon^{2(m+1)}}{n_0^{2r(m)}}
    +\alpha^2 C \Big\{\varepsilon^{2(m+1)} +\frac{1}{n} + \norm{\mathcal{Y}^{2,\varepsilon,(m,n)}_{T}}_2^2 \Big\},
\end{equation}
for all $\varepsilon,\alpha \in (0,1)$ and $n_0,n\geq 1$. 
\end{cor}

\end{rem}

\begin{rem}
If the driver of a BSDE contains
a linear part, we can transform the BSDE to the one considered in the current paper,
namely the equation (\ref{BSDE}).
For instance, let us solve FSDE (\ref{FSDE}) and the following BSDE:
\bea
-dY_t^{\varepsilon,\alpha} &=& [A(t,X_t^\varepsilon) Y_t^{\varepsilon,\alpha}  + Z_t^{\varepsilon}  B(t,X_t^\varepsilon) + 
\alpha f(t, X_t^\varepsilon,Y_t^{\varepsilon,\alpha},Z_t^{\varepsilon,\alpha})
]dt - Z_t^{\varepsilon} dW_t, \nn \\ 
Y_T^{\varepsilon,\alpha}  &=&g(X_T^\varepsilon),\nn
\eea
where 
$A$ and $B$ are $\mathbb{R}$-valued and $\mathbb{R}^d$-valued bounded functions,
respectively.

Let $\hat{Y}_t^{\varepsilon,\alpha}  := e^{\int_0^t A(s,X_t^\varepsilon) ds}Y_t^{\varepsilon,\alpha}$, $t\geq 0$,  
and $\hat{W}_t := W_t - \int_0^t B(s,X_s^\varepsilon) ds$, $t\geq 0$,
which is a Brownian motion under a probability measure $\hat{P}$
obtained by the change of measure with process $B(\cdot,X_{\cdot}^{\varepsilon})$.
Then, we have
{\footnotesize{\begin{align}
dX_t^\varepsilon
&=[\mu(t,X_t^\varepsilon)+
\varepsilon\sigma(t,X_t^\varepsilon)B(t,X_t^\varepsilon)]
dt+\varepsilon \sigma(t,X_t^\varepsilon)d\hat{W}_t,\quad
X_0^\varepsilon \in L^2(\Omega;\mathbb{R}^d),\nn\\
-d\hat{Y}_t^{\varepsilon,\alpha} &= \alpha e^{\int_0^t A(s,X_s^\varepsilon) ds} 
f(t, X_t^\varepsilon,Y_t^{\varepsilon,\alpha},Z_t^{\varepsilon,\alpha}) dt -  
e^{\int_0^t A(s,X_s^\varepsilon) ds} Z_t d\hat{W}_t
, \ Y_T^{\varepsilon,\alpha}  =  e^{\int_0^T A(s,X_s^\varepsilon) ds}g(X_T^\varepsilon).\nn
\end{align}}}
\end{rem}

\section{Numerical results}
In the numerical examples, we demonstrate 
that the deep BSDE method with the first order asymptotic expansion obtained in Proposition \ref{first_order_ae_U} provides enough accuracy in solving semilinear PDEs. The dimension $d$ in (\ref{BSDE}) is assumed to be $d=1$ or $d=100$. 

We investigate the accuracy of the new method by comparing to the standard Deep BSDE method in Weinan E et al. (2017) \cite{EHJ} and the Deep BSDE method with a prior knowledge in Fujii et al. (2019) \cite{FTT} for the model (\ref{BSDE}), where the target BSDEs with FSDEs are specified later. 

\subsection{Numerical schemes used in experiments}\label{sec-model}
In this subsection, we explain the details of schemes used in numerical experiments. To construct the deep neural networks for each method, we follow Weinan E et al. (2017) \cite{EHJ} and employ the adaptive moment estimation (Adam) with mini-batches. The parameters for the networks are set as follows: there are $d+10$ of hidden layers except batch normalization layers. For all learning steps, $256$ sample paths are generated and the learning rate is taken as  $0.01$.
\y

\noindent
{\bf (Numerical scheme)}
Now, let us briefly explain the schemes used in the numerical experiment in the following subsections.
\begin{enumerate}
\item 
{\bf (Deep BSDE method based on Weinan E et al. (2017))}
In forward discretization of $Y^{\varepsilon,\alpha}$, the Euler-Maruyama scheme $\bar{X}^{\varepsilon,(n)}$ is applied  with time step $n=20$. The initial guess of $Y_0^{\varepsilon,\alpha}$ is generated by uniform random number around $\overline{{\cal Y}}_0^{1,\varepsilon,(1)}$, which is a prior knowledge for the Deep BSDE method. 

In the study of Weinan E et al. (2017), it is known that the estimated value by the Deep BSDE method converges to the true value of $Y_0^{\varepsilon,\alpha}$ if we take a sufficient number of iteration steps. 

In Section \ref{sec-num-e} below, the estimate values based on this scheme are shown by the green lines labeled with ``Deep BSDE" in the figures. 

\item
{\bf (Deep BSDE method with an enhanced version of Fujii et al. (2019))}
In forward discretization of $Y^{\varepsilon,\alpha}$ in the Deep BSDE solver, as an approximation of
$\mathcal{Z}^{1,\varepsilon}$ we apply $\overline{{\cal Z}}_t^{1,\varepsilon,(1,n)}={\cal V}^{1,(1)}(t,\bar{X}_t^{\varepsilon,(n)})$, $t\geq 0$, with the function ${\cal V}^{1,(1)}$ defined by (\ref{cal V-1}) and the Euler-Maruyama scheme $\bar{X}^{\varepsilon,(n)}$ with time step $n=20$ to obtain an estimate of $\mathcal{Z}^{2,\varepsilon}$ by optimization in the Deep BSDE solver. As the initial value of $Y_0^{\varepsilon,\alpha}$, we use ${\calu}^{1,(1)}(0,x)$ with the function $\calu^{1,(1)}$ defined by (\ref{cal-U}), an approximation of $\mathcal{Y}^{1,\varepsilon}$, which appears in the linear part of our decomposition of the BSDE $(Y^{\varepsilon,\alpha},Z^{\varepsilon,\alpha})$ with $Y^{\varepsilon,\alpha}=\mathcal{Y}^{1,\varepsilon}+\alpha \mathcal{Y}^{2,\varepsilon}$ and $Z^{\varepsilon,\alpha}=\mathcal{Z}^{1,\varepsilon}+\alpha \mathcal{Z}^{2,\varepsilon}$. Thus, the scheme is an improved version of Fujii et al. (2019) \cite{FTT}, since it applies the higher order term $\overline{{\cal Z}}^{1,\varepsilon,(1)}$ than the leading order term $\overline{{\cal Z}}^{1,\varepsilon,(0)}$ that Fujii et al. (2019) \cite{FTT} uses. 

Through the study of Fujii et al. (2019) \cite{FTT}, it is also known that the estimated value by the enhanced Deep BSDE method converges to the true value of $Y_0^{\varepsilon,\alpha}$ with a much smaller number of iteration steps than by 
the original Deep BSDE method in Weinan E et al. (2017). 

In Section \ref{sec-num-e} below, the estimate values based on this scheme are shown by the red lines labeled with ``Deep BSDE[$(Y,Z)$]+AE[$\overline{\caly}^{1,\vep}_0$ and $\overline{Z}^{1,\vep}$]" in the figures.

\item
{\bf (New scheme)} Following the main result introduced in Section \ref{new-method}, particularly
Theorem \ref{main-thm}, we employ our approximation (\ref{new_decomp}) for the decomposition $Y^{\varepsilon,\alpha}_0=\mathcal{Y}^{1,\varepsilon}_0+\alpha \mathcal{Y}^{2,\varepsilon}_0$ with $m=1$ and $n=20$, namely,
\begin{eqnarray*}
Y_0^{\varepsilon,\alpha} &\approx& \overline{{\cal Y}}_0^{1,\varepsilon,(1)} + \alpha {\cal Y}_0^{2,\varepsilon,(1,20)\ast},\label{new-decomp2}
\end{eqnarray*}
where we compute the nonlinear part ${\cal Y}_0^{2,\varepsilon,(1,20)\ast}$ with 
(\ref{DBSDE_optimization_prob})--(\ref{linear-ae}) by the Deep BSDE solver,
while $\overline{{\cal Y}}_0^{1,\varepsilon,(1)}$
by ${\calu}^{1,(1)}(0,x)$ with the function $\calu^{1,(1)}$ defined by (\ref{cal-U}).

Specifically, in computation of ${\cal Y}_0^{2,\varepsilon,(1,20)\ast}$ by the Deep BSDE solver with the equation:
\beas
-d\mathcal{Y}_t^{2,\varepsilon}&= f(t,X_t^\varepsilon,\mathcal{Y}_t^{1,\varepsilon}+\alpha \mathcal{Y}_t^{2,\varepsilon},\mathcal{Z}_t^{1,\varepsilon}+ \alpha \mathcal{Z}_t^{2,\varepsilon})dt -  \mathcal{Z}_t^{2,\varepsilon} dW_t,\quad
\mathcal{Y}_T^{2,\varepsilon}=0, \label{Y2-2}
\eeas
we use $\overline{{\cal Y}}_t^{1,\varepsilon,(1)}={\cal U}^{1,(1)}(t,X_t^\varepsilon)$ and $\overline{{\cal Z}}_t^{1,\varepsilon,(1)}={\cal V}^{1,(1)}(t,X_t^\varepsilon)$ as approximations for $\mathcal{Y}_t^{1,\varepsilon}$ and $\mathcal{Z}_t^{1,\varepsilon}$ in the driver $f$, respectively.

In Section \ref{sec-num-e} below, the estimated values based on this new scheme  are shown by the blue lines labeled with ``New method [$\overline{{\cal Y}}^{1,\varepsilon}+{\cal Y}^{2,\varepsilon, DL}$, $\overline{{\cal Z}}^{1,\varepsilon}+{\cal Z}^{2,\varepsilon, DL}$]" in the figures. 
\end{enumerate}

The initial value of $Z_0^{\varepsilon,\alpha}$ or $\mathcal{Z}_0^{2,\varepsilon}$ is generated by uniform random number with the range $[-0.01,0.01]$ for each method. Numerical experiments presented in the following subsections are implemented by Python with TensorFlow on Google Colaboratory.

\subsection{Numerical experiments}\label{sec-num-e}

We show some numerical examples which show that our proposed method substantially outperforms other methods
in terms of terminal errors (numerical values of loss functions), variations and convergence speed. 

\subsubsection{The case of $d=1$}
This subsection presents the numerical results for the case of $d=1$.  \\

We first check the performance of our method in the model,
 where the explicit value of the solution is obtained by the Picard iteration. We consider an option pricing model in finance that takes CVA(credit value adjustment) into account as follows: 
\begin{align}
  dX_t^\varepsilon =&\mu(t,X_t^\varepsilon)dt+\varepsilon \sigma(t,X_t^\varepsilon)dW_t,\\
  -dY_t^{\varepsilon,\alpha}
  =&-\alpha (Y_t^{\varepsilon,\alpha})^+ dt - Z_t^{\varepsilon,\alpha} dW_t,\quad
  Y_T^{\varepsilon,\alpha}=(X_T^\varepsilon-K)^+,
\end{align}
with
$f(t,x,y,z)=-(y)^+$ and $g(x)=(x-K)^+$.
in (\ref{FSDE}) and (\ref{BSDE}).
We note that $\alpha=$(loss rate in default)$\times$(default intensity)
in a finance model of CVA. 

In computation we set 
  $\mu(t,x)=0$, $\sigma(t,x)=x$, $\varepsilon=\sigma=0.2$, $X_0 = 100$,
   $\alpha=0.05$, 
  $T=0.5$ with $K=100$ (ATM case) and $K=115$ (OTM case).
  
 In this case an explicit value of $Y_0$ is computed as $\textstyle{Y_0={\cal Y}_0^1(1+\sum_{i=1}^\infty (-1)^i \alpha^i T^i \frac{1}{i!})}$. More precisely, by the $k$-Picard iteration of the backward equation: 
\begin{align}
&
 -dY_t^{\varepsilon,\alpha,[k]}
  =-\alpha (Y_t^{\varepsilon,\alpha,[k-1]})^+ dt - Z_t^{\varepsilon,\alpha,[k]} dW_t, \ \ \ Y_T^{\varepsilon,\alpha,[k]}=(X_T^\varepsilon-K)^+,\\
\mbox{with} & \ \ \ \
  (Y_t^{\varepsilon,\alpha,[0]})^+=E[(X_T^\varepsilon-K)^+|{\cal F}_t]={\cal Y}_t^1, \ \mbox{for \ all} \ t\geq 0, 
\end{align}
it is easy to see that $\textstyle{Y_t^{\varepsilon,\alpha,[k]}={\cal Y}_t^1-\alpha \int_t^TE [Y_s^{\varepsilon,\alpha,[k-1]} |{\cal F}_t ]ds}$, and thus one has
\beas
 \textstyle{Y_0^{\varepsilon,\alpha,[k]}={\cal Y}_0^1(1+\sum_{i=1}^k (-1)^i \alpha^i T^i \frac{1}{i!})}.  
\eeas

Then, the true values are given as $Y_0=5.50$ in the ATM case and $Y_0=1.26$ in the OTM case by the $5$-Picard iteration, which provides enough convergence 
and hence accuracy. 

Figure 1 and 2 show the numerical values of loss functions and the approximate values of $Y_0$, respectively against the number of iteration steps for the ATM case, and Figure 3 and 4 for the OTM case. 

By Figure 2 and 4, the numerical values of ``New method [$\overline{{\cal Y}}^{1,\varepsilon}+{\cal Y}^{2,\varepsilon, DL}$,$\overline{{\cal Z}}^{1,\varepsilon}+{\cal Z}^{2,\varepsilon, DL}]$" converge to the true values substantially faster with smaller variations comparing to other schemes. Also, we can see that the errors of ``New method" are much smaller according to the behavior of their loss functions against the number of iteration steps in Figure 1 and 3.    

\setcounter{figure}{0}
\begin{figure}[H]
  \centering
  \includegraphics[height=8.5cm]{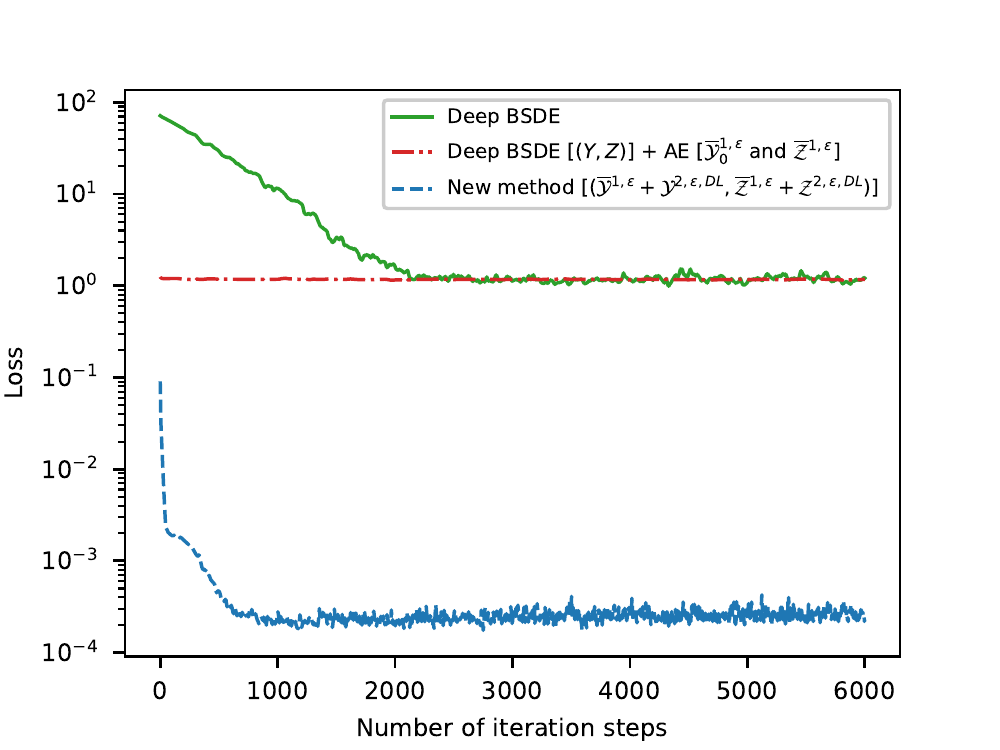}
  \caption{Values of the loss function and number of iteration steps (1-dim option pricing model with CVA, ATM case)}
\end{figure}

\begin{figure}[H]
  \centering
  \includegraphics[height=8.5cm]{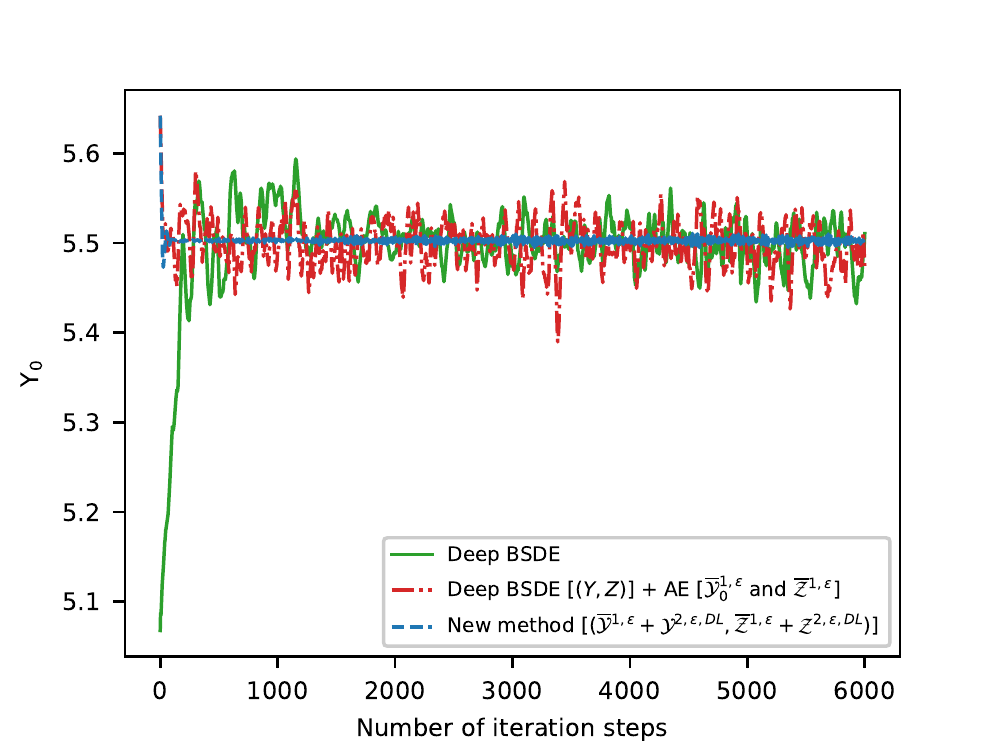}
  \caption{ Approximate values of $Y_0$ (true value: 5.50) and number of iteration steps (1-dim option pricing model with CVA, ATM case)}
\end{figure}

\begin{figure}[H]
  \centering
  \includegraphics[height=8.5cm]{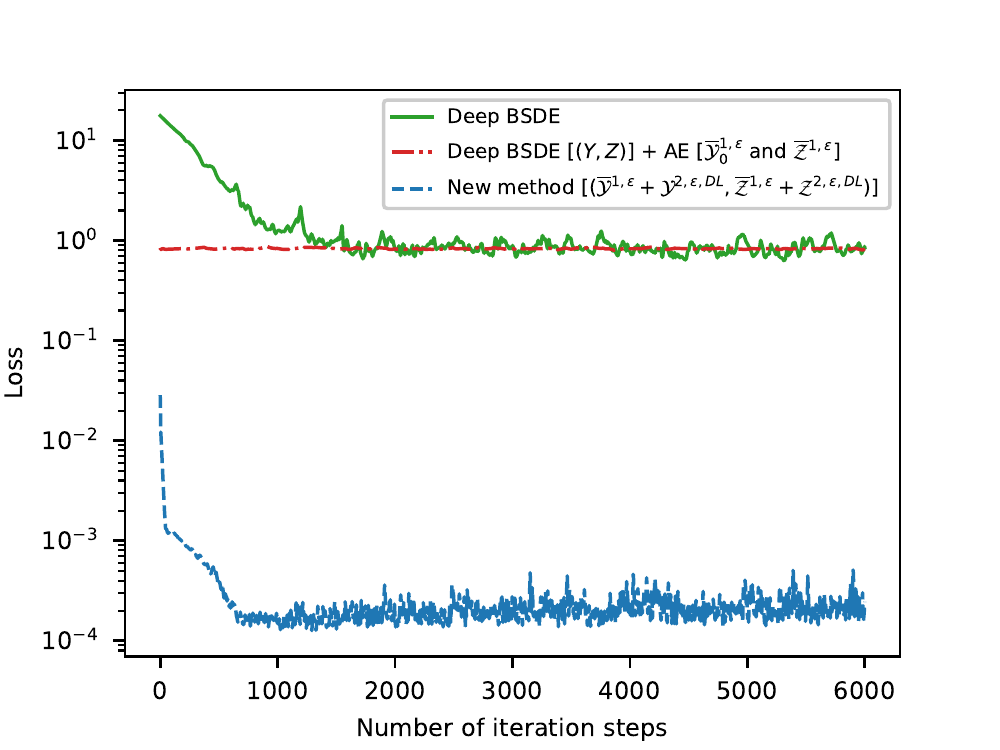}
  \caption{Values of the loss function and number of iteration steps (1-dim option pricing model with CVA, OTM case)}
\end{figure}

\begin{figure}[H]
  \centering
  \includegraphics[height=8.5cm]{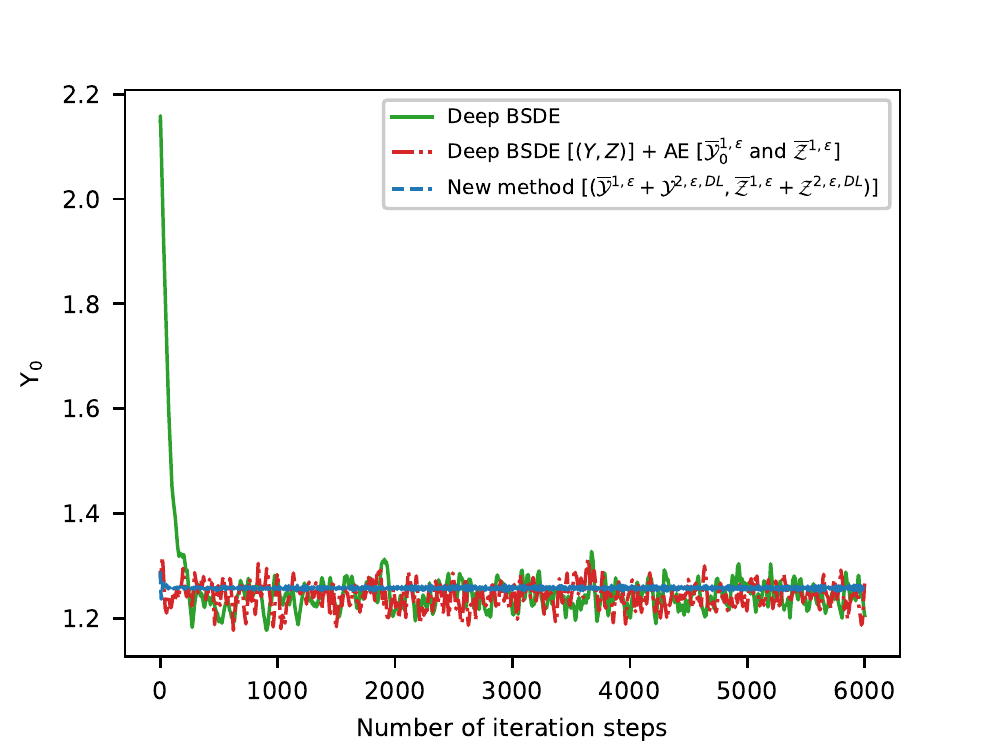}
  \caption{ Approximate values of $Y_0$ (true value: 1.26) and number of iteration steps (1-dim option pricing model with CVA, OTM case)}
\end{figure}

Next, we present numerical examples for the model, where explicit values of $Y_0$ can not be obtained without numerical schemes such as
Monte Carlo simulations. Let us consider 
\begin{align}
  dX_t^\varepsilon =&\mu(t,X_t^\varepsilon)dt+\varepsilon \sigma(t,X_t^\varepsilon)dW_t,\\
  -dY_t^{\varepsilon,\alpha}
  =&\alpha f(t,X_t^\varepsilon,Y_t^{\varepsilon,\alpha},Z_t^{\varepsilon,\alpha}) dt - Z_t^{\varepsilon,\alpha} dW_t,\quad
  Y_T^{\varepsilon,\alpha}=g(X_T^\varepsilon)
\end{align}
with $\mu(t,x)=0$, $\sigma(t,x)=x$, $\varepsilon=\sigma=0.2$, $T=0.25$,
$X_0=100$,
and
\begin{align}
  f(t,x,y,z)
  &= -\left\{ \left(y - z \sigma^{-1} \mathrm{1} \right)^- (R-r)\right\}
\end{align}
with $R=0.06$, $r=0.0$, and 
\begin{equation}
	g(x)=(x - K_1)^+- 2 (x - K_2)^+, 
	\ \ \mbox{with} \ K_1=95, \ K_2=105.
\end{equation}

As we mentioned in Section \ref{sec-model}, the estimated value by the methods ``Deep BSDE" and ``Deep BSDE[$(Y,Z)$]+AE[$\overline{\caly}^{1,\vep}_0$ and $\overline{Z}^{1,\vep}$]" converges to the true value of $Y_0$. Then, in the experiments, we check whether the estimated value by ``New method [$\overline{{\cal Y}}^{1,\varepsilon}+{\cal Y}^{2,\varepsilon, DL}$,$\overline{{\cal Z}}^{1,\varepsilon}+{\cal Z}^{2,\varepsilon, DL}]$" converges faster than the ones computed by the methods ``Deep BSDE" and ``Deep BSDE[$(Y,Z)$]+AE[$\overline{\caly}^{1,\vep}_0$ and $\overline{Z}^{1,\vep}$]". 

Figure 5 shows the numerical values of loss functions against the number of iteration steps. 
While ``Deep BSDE[$(Y,Z)$]+AE[$\overline{\caly}^{1,\vep}_0$ and $\overline{Z}^{1,\vep}$]" is superior to the original ``Deep BSDE", we see that ``New method [$\overline{{\cal Y}}^{1,\varepsilon}+{\cal Y}^{2,\varepsilon, DL}$,$\overline{{\cal Z}}^{1,\varepsilon}+{\cal Z}^{2,\varepsilon, DL}]$" gives much more stable and accurate convergence than other schemes. 

Figure 6 shows the approximate values of $Y_0$ against the number of iteration steps. It is observed that ``New method [$\overline{{\cal Y}}^{1,\varepsilon}+{\cal Y}^{2,\varepsilon, DL}$,$\overline{{\cal Z}}^{1,\varepsilon}+{\cal Z}^{2,\varepsilon, DL}]$" provides the fastest convergence with the smallest standard deviation, while ``Deep BSDE[$(Y,Z)$]+AE[$\overline{\caly}^{1,\vep}_0$ and $\overline{Z}^{1,\vep}$]" gives better approximation than ``Deep BSDE". 

\begin{figure}[H]
  \centering
  \includegraphics[height=8.5cm]{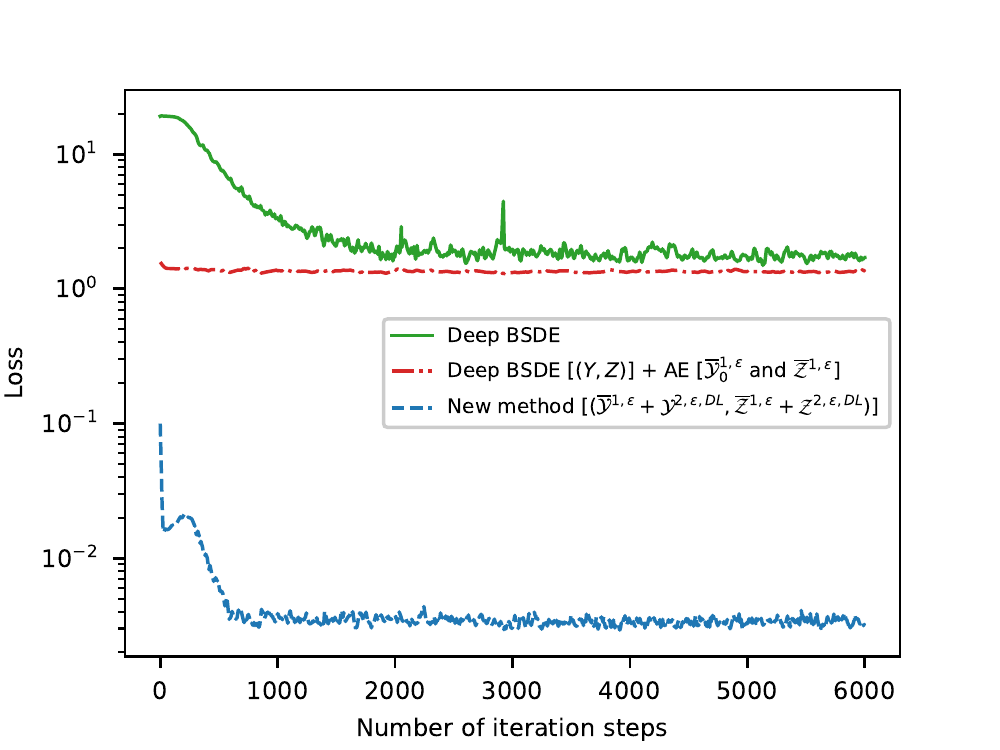}
  \caption{Values of the loss function and number of iteration steps}
\end{figure}

\begin{figure}[H]
  \centering
  \includegraphics[height=8.5cm]{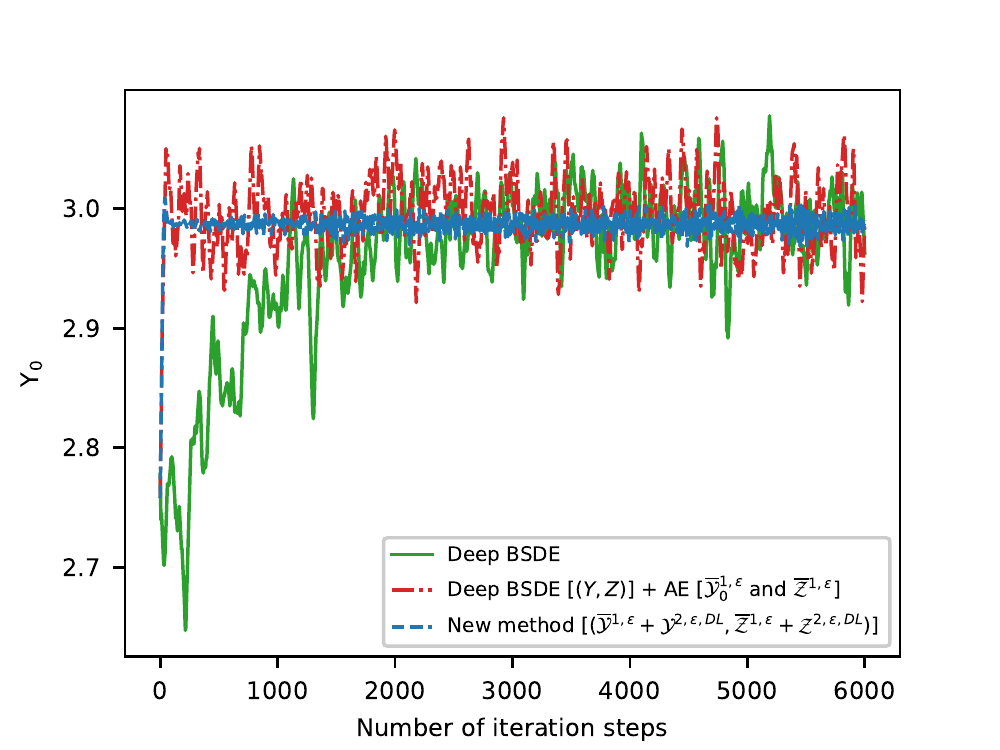}
  \caption{ Approximate values of $Y_0$ and number of iteration steps}
\end{figure}

\subsubsection{The case of $d=100$}
We show the main numerical result for $d=100$. The same experiment as in the case of $d=1$ is performed. Let us consider
\begin{align}
  dX_t^{\varepsilon,i} =&\mu^i(t,X_t^\varepsilon)dt+\varepsilon \sum_{j=1}^d \sigma^i_j(t,X_t^\varepsilon)dW^j_t,\\
  -dY_t^{\varepsilon,\alpha}
  =&\alpha f(t,X_t^\varepsilon,Y_t^{\varepsilon,\alpha},Z_t^{\varepsilon,\alpha}) dt - Z_t^{\varepsilon,\alpha} dW_t,\quad
  Y_T^{\varepsilon,\alpha}=g(X_T^\varepsilon),
\end{align}
with $\mu^i(t,x)=0$, $\sigma^i_j(t,x)=x^i$ ($i=1,\cdots,d$), $\varepsilon=0.4$,
$X_0^i=100$, $T=0.25$, and
\begin{align}
  f(t,x,y,z)
  &= -\left\{ 
   - \left(y - \sum_{k=1}^d \sum_{j=1}^d z_k [\sigma^{-1}]_{kj} \right)^- (R-r)\right\},
\end{align}
with $R=0.01$, $r=0.0$ where $\theta$ is defined by $\mu-r\mathrm{1}=\sigma \theta$, and 
\begin{equation}
  g(x)
  =\left(\frac{1}{d} \sum_{i=1}^d x_i - K_1\right)^+
  - 2 \left(\frac{1}{d} \sum_{i=1}^d x_i - K_2\right)^+, \ \ \mbox{with} \ K_1=95, \ K_2=105.
\end{equation} 
The result is given in Figure 7, 8 and 9. 
 It seems that the convergence speed of the original deep BSDE method is too slow to obtain the precise result. On the contrary, ``Deep BSDE[$(Y,Z)$]+AE[$\overline{\caly}^{1,\vep}_0$ and $\bar{Z}^{1,\vep}$]" and ``New method [$\overline{{\cal Y}}^{1,\varepsilon}+{\cal Y}^{2,\varepsilon, DL}$,$\overline{{\cal Z}}^{1,\varepsilon}+{\cal Z}^{2,\varepsilon, DL}]$" work well even in this high dimensional case. Particularly, 
``New method [$\overline{{\cal Y}}^{1,\varepsilon}+{\cal Y}^{2,\varepsilon, DL}$,$\overline{{\cal Z}}^{1,\varepsilon}+{\cal Z}^{2,\varepsilon, DL}]$" provides a remarkable performance in terms of convergence speed, accuracy (numerical values of loss functions) and variations. 
Moreover, comparing the results of our new method and
``Deep BSDE[$(Y,Z)$]+AE[$\overline{\caly}^{1,\vep}_0$ and $\bar{Z}^{1,\vep}$]" closely,
Figure 9 shows that the variation of $Y_0$ by our method
is much smaller,
which is consistent with much smaller values of loss functions for the new method appearing in Figure 7.

\begin{figure}[H]
  \centering
  \includegraphics[height=8.5cm]{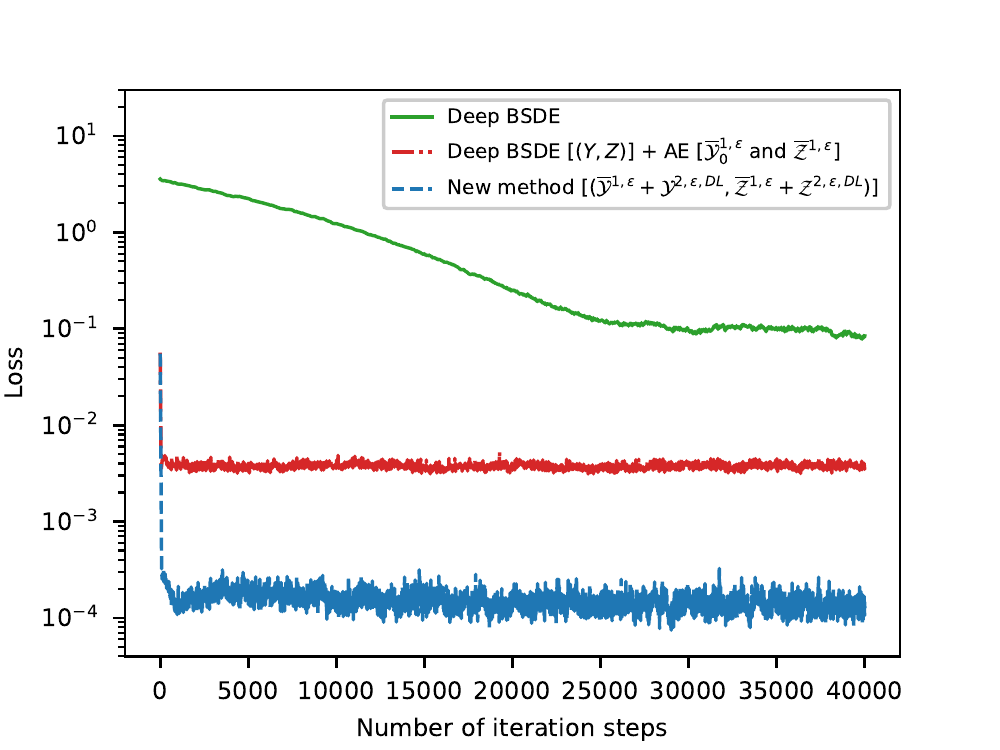}
  \caption{Values of the loss function and number of iteration steps}
\end{figure}

\begin{figure}[H]
  \centering
  \includegraphics[height=8.5cm]{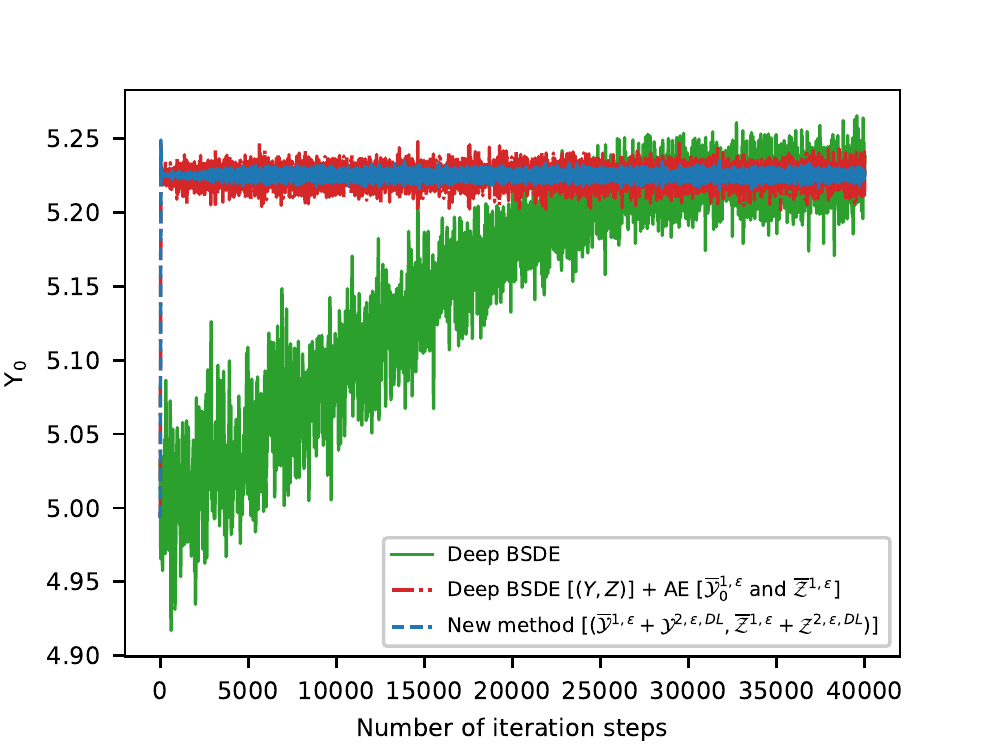}
  \caption{Approximate values of $Y_0$ and number of iteration steps}
\end{figure}

\begin{figure}[H]
  \centering
  \includegraphics[height=8.5cm]{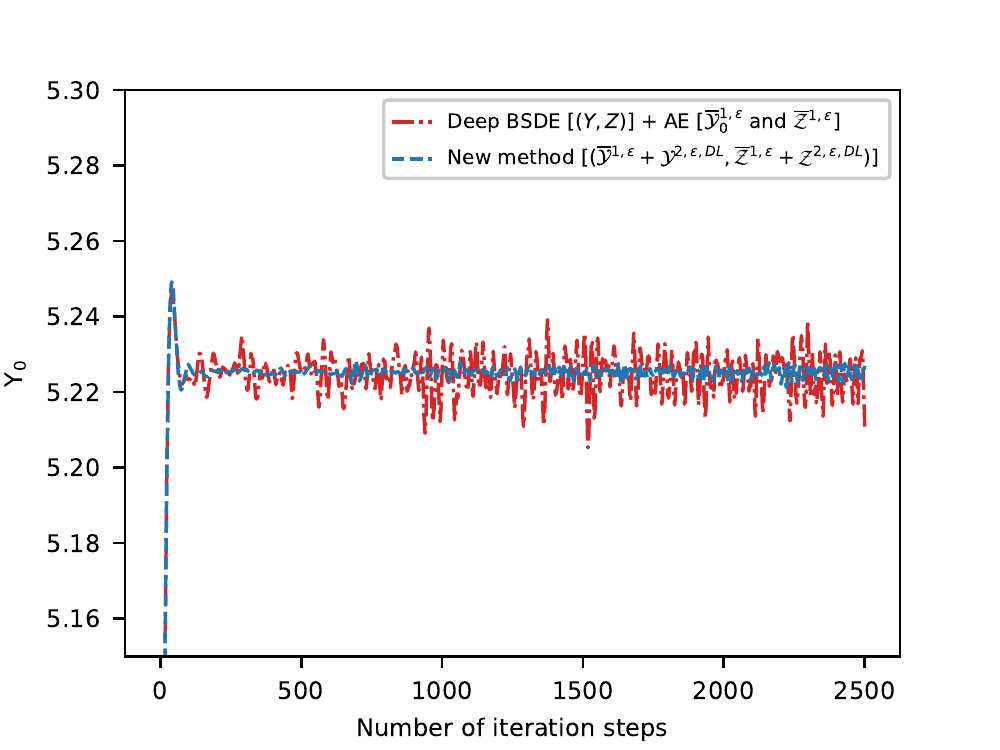}
  \caption{Approximate values of $Y_0$ and number of iteration steps (enlarged view for ``Deep BSDE[$(Y,Z)$]+AE[$\overline{\caly}^{1,\vep}_0$ and $\overline{Z}^{1,\vep}$]'' and ``New method [$\overline{{\cal Y}}^{1,\varepsilon}+{\cal Y}^{2,\varepsilon, DL}$, $\overline{{\cal Z}}^{1,\varepsilon}+{\cal Z}^{2,\varepsilon, DL}$]'')}
\end{figure}

\section{Conclusion and future works}
This paper has introduced a new control variate method for Deep BSDE solver to improve the methods such as in Weinan E et al. (2017) \cite{EHJ} and Fujii et al. (2019) \cite{FTT}. First, we decompose a target semilinear PDE (BSDE) to two parts, namely 
dominant linear and residual nonlinear PDEs (BSDEs). 
When the dominant part is obtained as a closed-form or approximated based on an asymptotic expansion scheme, the small nonlinear PDE part is efficiently computed by Deep BSDE solver, where the asymptotic expansion crucially works as a control variate. 
 The main theorem provides the validity of our proposed method. Moreover, numerical examples for one and 100 dimensional BSDEs 
corresponding to target nonlinear PDEs show the effectiveness of our scheme, which is consistent with our initial conjecture and theoretical result.

As mentioned in Remark \ref{remark_1}, 
even if the accuracy of the standard asymptotic expansion scheme becomes worse,
the linear PDE part can be more efficiently approximated by the existing methods such as \cite{TYo}\cite{TY16}\cite{Y19}\cite{NYmcma}\cite{OY}\cite{IY}. We should check those performances in such cases against various nonlinear models. Also, it will be a challenging task to examine whether the high order automatic differentiation schemes proposed in \cite{YY19}\cite{ToYa} work as efficient approximations of $Z$ in nonlinear BSDEs or $\partial_x u$ in nonlinear PDEs. These are left for future studies. 

\appendix
\section*{Appendix}
\section{Proof of Propositions}
\subsection{Proof of Proposition \ref{Prop_expansion_estimate}}\label{Proof_Prop_expansion_estimate}
See Proposition 4.2 in Takahashi and Yamada (2015) \cite{TY15} for (\ref{ae_error_1})--(\ref{ae_error_2-2}), for instance. 

Also, note that for $p\geq 1$ and a multi-index $\alpha$, $\sup {}_{x\in \mathbb{R}^d}\|\partial_x^{\alpha} \overline{X}^{t,x,\varepsilon}_T\|_p \leq C(T)$ and $\sup {}_{x\in \mathbb{R}^d} \|\partial_x^{\alpha} {\cal W}^{t,x,\varepsilon,(m)}_{T}\|_p \leq C(T)$ hold for $t<T$. Then, $\sup {}_{x\in \mathbb{R}^d} |\partial_x^{2} {\cal U}^{1,(m)}(t,x)| \leq \| \nabla^2 g \|_{\infty} C(T)$ for $t<T$,
i.e. ${\cal U}^{1,(m)}(t,\cdot) \in C_b^2(\mathbb{R}^d)$, $t<T$.

 For ${\cal V}^{1,(m)}$, we have the representation 
\begin{align*}
{\cal V}^{1,(m)}(t,x)=E[ g(\overline{X}^{t,x,\varepsilon}_T)  {\cal Z}^{t,x,\varepsilon,(m)}_{T} ]=E[ (\nabla g)(\overline{X}^{t,x,\varepsilon}_T)  {\cal Q}^{t,x,\varepsilon,(m)}_{T} ],
\end{align*}
 for a matrix-valued Wiener functional ${\cal Q}^{t,x,\varepsilon,(m)}_{T}=[ [{\cal Q}^{t,x,\varepsilon,(m)}_{T}]^i_j ]_{1\leq i,j \leq d}$ such that $[{\cal Q}^{t,x,\varepsilon,(m)}_{T}]^i_j \in \mathbb{D}^\infty$, $1\leq i,j \leq d$, satisfying for $p\geq 1$ and a multi-index $\alpha$, $\sup {}_{x\in \mathbb{R}^d}\|\partial_x^{\alpha} {\cal Q}^{t,x,\varepsilon,(m)}_{T}\|_p \leq C(T)$ for $t<T$. Then, $\sup {}_{x\in \mathbb{R}^d} |\partial_x {\cal V}^{1,(m)}(t,x)| \leq \| \nabla^2 g \|_{\infty} C(T)$ for $t<T$, i.e.
 ${\cal V}^{1,(m)}(t,\cdot) \in C_b^1(\mathbb{R}^d)$, $t<T$.
 $\Box$ 

\subsection{Proof of Proposition \ref{first_order_ae_U}}\label{Proof_Prop_expansion_order1}
For the derivations, we use Malliavin calculus. Let ${\cal T} \in \mathcal{S}'(\mathbb{R}^d)$ be a tempered distribution and $F \in (\mathbb{D}^\infty)^d$ be a nondegenerate Wiener functional in the sense of Malliavin. 
Then, ${\cal T} (F)$ is well-defined as an element of the space of Watanabe distributions $\mathbb{D}^{-\infty}$, that is the dual space of $\mathbb{D}^{\infty}$. Also, for $G \in \mathbb{D}^{\infty}$, a (generalized)
expectation $E[{\cal T}(F)G]$ is understood as 
a coupling of ${\cal T}(F)\in \mathbb{D}^{-\infty}$
and $G\in \mathbb{D}^{\infty}$, namely
${}_{\mathbb{D}^{-\infty}}\langle {\cal T}(F),G \rangle_{\mathbb{D}^{\infty}}$. \\

Note that $G_T^{t,x,\varepsilon}:=(X_T^{t,x,\varepsilon}-X_T^{t,x,0})/\varepsilon$ and $(\partial/\partial x)X_T^{t,x,\varepsilon}$ in
\begin{align}
{\cal U}^1(t,x)=E[g(X_T^{t,x,\varepsilon})] =\int_{\mathbb{R}^d}g(X_T^{t,x,0}+\varepsilon y) E[ \delta_y(G_T^{t,x,\varepsilon}) ]dy 
\label{integral1}
\end{align}
 and 
\begin{align} 
& (\partial/\partial x){\cal U}^1\sigma^{\varepsilon} (t,x) 
 = E[(\nabla g)(X_T^{t,x,\varepsilon})(\partial/\partial x)X_T^{t,x,\varepsilon}] \varepsilon \sigma(t,x) \\
 &=\int_{\mathbb{R}^d} \sum_{i=1}^d (\partial_i g)(X_T^{t,x,0}+\varepsilon y) E[ \delta_y(G_T^{t,x,\varepsilon}) (\partial/\partial x)X_T^{t,x,\varepsilon,i} ] dy \varepsilon \sigma(t,x)
  \label{integral2}
\end{align} 
 have expansions in $\mathbb{D}^\infty$ whose expansion coefficients are given by iterated stochastic integrals: 
 $G_T^{t,x,\varepsilon} \sim X_{1,T}^{t,x}+ \varepsilon X_{2,T}^{t,x}+\cdots$ and $(\partial/\partial x)X_T^{t,x,\varepsilon} \sim J_{t \to T}^{0,x} +\varepsilon J_{t \to T}^{1,x}+ \cdots$. In particular,  
  \begin{align}
    X_{1,T}^{t,x}
    &= \sum_{k=1}^d \int_t^T J_{t \to T}^{0,x}(J_{t \to s}^{0,x})^{-1} \sigma_k(X_s^{t,x,0}) dW_s^k,\\
    X_{2,T}^{t,x}
    &=\sum_{k=1}^d \int_t^T J_{t \to T}^{0,x}(J_{t \to s}^{0,x})^{-1} \partial \sigma_k(X_s^{t,x,0}) X_{1,s}^{t,x} dW_s^k\\
    & \ \ \ \ +\frac{1}{2} \int_t^T J_{t \to T}^{0,x}(J_{t \to s}^{0,x})^{-1} \partial^2 \mu(X_s^{t,x,0}) \cdot X_{1,s}^{t,x} \otimes X_{1,s}^{t,x} ds,
  \end{align} 
  and 
  \begin{align}
  &J_{t \to T}^{0,x}=(\partial/\partial x)X_T^{t,x,0},\\
    &J_{t \to T}^{1,x}
    =J_{t \to T}^{0,x}
    \{
      \int_t^T \partial^2 \mu(X_s^{t,x,0}) X_{1,s}^{t,x} ds
      + \sum_{k=1}^d \int_t^T \partial \sigma_k(X_s^{t,x,0}) dW_s^k
    \},
  \end{align}
  where the followings are used: for a smooth function $V:\mathbb{R}^d \to \mathbb{R}^d$,
  \begin{equation}
    \partial^2 V(x)
    =\left[\pdv[2]{V_\alpha^i(x)}{x^j}{x^k}\right]^k_j,
  \end{equation}
  \begin{equation}
    \left[\partial^2 V \cdot \eta \otimes \eta\right]^i
    =\sum_{j,k} \pdv[2]{V_\alpha^i}{x^j}{x^k} \eta^j \eta^k, \ \ \ \eta \in \mathbb{R}^d.
  \end{equation}
 
We expand $E[ \delta_y(X_T^{t,x,\varepsilon}) ]$ in (\ref{integral1}) and $E[ (\nabla \delta_y)(X_T^{t,x,\varepsilon}) (\partial/\partial x)X_T^{t,x,\varepsilon} ]$ in (\ref{integral2}) to obtain explicit expressions of ${\cal U}^{1,(1)}(t,x)$ and ${\cal V}^{1,(1)}(t,x)$. 
Next, let us recall the following formulas. 

\begin{lemma}\label{lemma_appendix}
Let ${\cal T} \in \mathcal{S}'(\mathbb{R}^d)$ be a tempered distribution. 
\noindent
\begin{enumerate}
\item 
 For an adapted process $h \in L^2([0,T] \times \Omega)$, 
\begin{align}
\sum_{j=1}^d E[ \partial_j {\cal T} ({X}_{1,T}^{t,x}) \int_t^T (D_{i,s}{X}_{1,T}^{t,x,j}) h(s) ds  ] = E[ {\cal T} ({X}_{1,T}^{t,x}) \int_t^T h(s)dW_s^i ],\label{mall_IBP}
\end{align}
where $D_{i,\cdot} F$ represents the $i$-th element of the Malliavin derivative \\ $D_{\cdot}F=(D_{1,\cdot}F,\cdots,D_{d,\cdot}F )$ for $F \in \mathbb{D}^\infty$. 
\item For $1 \leq i_1,\cdots,i_\ell \leq d$,
\begin{align}
E[ (\partial_{{i_1}} \cdots \partial_{{i_\ell}} {\cal T}) ({X}_{1,T}^{t,x}) ]=E[ {\cal T}({X}_{1,T}^{t,x}) H_{(i_1,\cdots,i_\ell)}({X}_{1,T}^{t,x},1) ]. \ \ \ \ \ \label{diff_computation}
\end{align}
\end{enumerate}
\end{lemma}
\noindent\\
{\it Proof of Lemma \ref{lemma_appendix}}. Use the duality formula (see Theorem 1.26 of \cite{MallThal} or Proposition 1.3.11 of \cite{N}), with $\textstyle{D{\cal T}(\Xi)=\sum_{i=1}^d (\partial_i {\cal T})(\Xi)D\Xi^i}$ for $\Xi=(\Xi^1,\cdots,\Xi^d) \in (\mathbb{D}^\infty)^d$ (see Proof of Proposition 2.1.9 of \cite{N} or Proof of Theorem 2.6 of  \cite{TY12}) to obtain the first assertion. Also, the second assertion is immediately obtained by the integration by parts. $\Box$

In the expansions of (\ref{integral1}) and (\ref{integral2}), iterated integrals such as 
{\footnotesize{\begin{align}
\int_t^T  h_{j_1}(t_1) \int_t^{t_1} h_{j_2}(t_2)dW_{t_2}^{j_2}dW_{t_1}^{j_1} \ \ (h_{j_l} \in L^2([0,T]), \ l=1,2, \ j_1,j_2=1,\cdots,d)
\end{align}}} 
appear, for which the following calculation holds with use of (\ref{mall_IBP}): 
{\footnotesize{\begin{align}
& \sum_{i_1} E[ \partial_{i_1} {\cal T} ({X}_{1,T}^{t,x}) \int_t^T  h_{j_1}(t_1) \int_t^{t_1} h_{j_2}(t_2)dW_{t_2}^{j_2}dW_{t_1}^{j_1} ] \label{calc_ibp} \\
=&\sum_{i_1,i_2} E[ \partial_{i_2} \partial_{i_1}{\cal T} ({X}_{1,T}^{t,x}) \int_t^T  (D_{j_1,t_1}{X}_{1,T}^{t,x,i_2})  h_{j_1}(t_1) \int_t^{t_1} h_{j_2}(t_2)dW_{t_2}^{j_2}dt_1 ]\nn\\
=&\sum_{i_1,i_2} \int_t^T  (D_{j_1,t_1}{X}_{1,T}^{t,x,i_2})  h_{j_1}(t_1) E[ \partial_{i_2}  \partial_{i_1} {\cal T} ({X}_{1,T}^{t,x})  \int_t^{t_1} h_{j_2}(t_2)dW_{t_2}^{j_2} ] dt_1 \nn\\
=&\sum_{i_1,i_2,i_3} \int_t^T  (D_{j_1,t_1}{X}_{1,T}^{t,x,i_2})  h_{j_1}(t_1) E[ \partial_{i_3} \partial_{i_2} \partial_{i_1}{\cal T} ({X}_{1,T}^{t,x})  \int_t^{t_1} (D_{j_2,t_2}{X}_{1,T}^{t,x,i_3})  h_{i_2}(t_2)d{t_2} ] dt_1\nn\\
=&\sum_{i_1,i_2,i_3}E[ \partial_{i_3} \partial_{i_2}\partial_{i_1} {\cal T} ({X}_{1,T}^{t,x}) ] \int_t^T  (D_{j_1,t_1}{X}_{1,T}^{t,x,i_2})  h_{j_1}(t_1) \int_t^{t_1} (D_{j_2,t_2}{X}_{1,T}^{t,x,i_3})  h_{j_2}(t_2)d{t_2} dt_1.\nn
\end{align}}}
Note that 
$s \mapsto D_{j,s}{X}_{1,T}^{t,x,i}$ is deterministic, and one has 
\begin{align}
D_{j,s}{X}_{1,T}^{t,x,i}=[J^{0,x}_{t \to T}{J^{0,x}_{t \to s}}^{-1}\sigma_j(s,X_s^{t,x,0})]^i.
\end{align}
 Thus, we get 
{\footnotesize{\begin{align}
&\sum_{i_1} E[ \partial_{i_1} {\cal T} ({X}_{1,T}^{t,x}) \int_t^T  h_{i_1}(t_1) \int_t^{t_1} h_{i_2}(t_2)dW_{t_2}^{i_2}dW_{t_1}^{i_1} ]\\
=&
 \sum_{i_1,i_2,i_3}E[ \partial_{i_3} \partial_{i_2} \partial_{i_1} {\cal T} ({X}_{1,T}^{t,x}) ] \\
&  \int_t^T  [J^{0,x}_{t \to T}(J^{0,x}_{t \to t_1})^{-1}\sigma_{j_1}(t_1,X_{t_1}^{t,x,0})]^{i_2}  h_{j_1}(t_1) \int_t^{t_1} [J^{0,x}_{t \to T}(J^{0,x}_{t \to t_2})^{-1}\sigma_{j_2}(t_2,X_{t_2}^{t,x,0})]^{i_3} h_{j_2}(t_2)d{t_2} dt_1.
\end{align}}} 
Using the above calculation with (\ref{diff_computation}), we have 
\begin{align}
&\sum_{i_1} E[ \partial_{i_1} {\cal T} ({X}_{1,T}^{t,x}) \varepsilon X_{2,T}^{t,x,i_1} ]\\
=&
\varepsilon \sum_{i_1,i_2,i_3,j_1,k_1,k_2} E[ \partial_{i_3}\partial_{i_2}\partial_{i_1} {\cal T} ({X}_{1,T}^{t,x})   ] C_{i_1,i_2,i_3,j_1}^{(1),k_1,k_2}(t,T,x)\\
&+\varepsilon \sum_{i_1,i_2,i_3,j_1,j_2,k_1,k_2} E[ \partial_{i_3}\partial_{i_2}\partial_{i_1} {\cal T} ({X}_{1,T}^{t,x})   ] C_{i_1,i_2,i_3,j_1,j_2}^{(2),k_1,k_2}(t,T,x)\\
&+\varepsilon \sum_{i_1,j_1,j_2,k_1,k_2} E[ \partial_{i_1} {\cal T} ({X}_{1,T}^{t,x})   ] \frac{1}{2} 1_{k_1=k_2} C_{i_1,j_1,j_2}^{(3),k_1,k_2}(t,T,x)\\
=&
\varepsilon \sum_{i_1,i_2,i_3,j_1,k_1,k_2} E[ {\cal T} ({X}_{1,T}^{t,x}) H_{(i_1,i_2,i_3)}({X}_{1,T}^{t,x},1)  ] C_{i_1,i_2,i_3,j_1}^{(1),k_1,k_2}(t,T,x)\\
&+\varepsilon \sum_{i_1,i_2,i_3,j_1,j_2,k_1,k_2} E[ {\cal T} ({X}_{1,T}^{t,x}) H_{(i_1,i_2,i_3)}({X}_{1,T}^{t,x},1)  ] C_{i_1,i_2,i_3,j_1,j_2}^{(2),k_1,k_2}(t,T,x)\\
&+\varepsilon \sum_{i_1,,j_1,j_2,k_1,k_2} E[ {\cal T} ({X}_{1,T}^{t,x}) H_{(i_1)}({X}_{1,T}^{t,x},1)  ] \frac{1}{2} 1_{k_1=k_2} C_{i_1,j_1,j_2}^{(3),k_1,k_2}(t,T,x).
\end{align}
Therefore, we get (\ref{cal-U}) as: 
\begin{align}
{\cal U}^{1,(1)}(t,x)
=&E[g(\overline{X}_T^{t,x,\varepsilon})]\\
&+\varepsilon \sum_{i_1,i_2,i_3,j_1,k_1,k_2} E[ g(\overline{X}_T^{t,x,\varepsilon}) H_{(i_1,i_2,i_3)}({X}_{1,T}^{t,x},1)  ] C_{i_1,i_2,i_3,j_1}^{(1),k_1,k_2}(t,T,x)\\
&+\varepsilon \sum_{i_1,i_2,i_3,j_1,j_2,k_1,k_2} E[ g(\overline{X}_T^{t,x,\varepsilon}) H_{(i_1,i_2,i_3)}({X}_{1,T}^{t,x},1)  ] C_{i_1,i_2,i_3,j_1,j_2}^{(2),k_1,k_2}(t,T,x)\\
&+\varepsilon \sum_{i_1,,j_1,j_2,k_1,k_2} E[g(\overline{X}_T^{t,x,\varepsilon}) H_{(i_1)}({X}_{1,T}^{t,x},1)  ] \frac{1}{2} 1_{k_1=k_2} C_{i_1,j_1,j_2}^{(3),k_1,k_2}(t,T,x).
\end{align}

Next, we give the representation (\ref{cal V-1}). The function $(\partial/\partial x){\cal U}^1\sigma^{\varepsilon}$ given by
\begin{align} 
& (\partial/\partial x){\cal U}^1\sigma^{\varepsilon} (t,x) 
 = E[(\nabla g)(X_T^{t,x,\varepsilon})(\partial/\partial x)X_T^{t,x,\varepsilon}] \varepsilon \sigma(t,x) \\
 &=\int_{\mathbb{R}^d} \sum_{i=1}^d (\partial_i g)(X_T^{t,x,0}+\varepsilon y) E[ \delta_y(G_T^{t,x,\varepsilon}) (\partial/\partial x)X_T^{t,x,\varepsilon,i} ] dy \varepsilon \sigma(t,x)
\end{align} 
is expanded as 
  \begin{align}
    &{\cal V}^{1,(1)}(t,x)\\
    &=\frac{1}{\varepsilon} E[ \sum_{i_1} g(X_T^{t,x,0} + \varepsilon X_{1,T}^{t,x}) H_{(i_1)}( X_{1,T}^{t,x},1 )
    ]J_{t \to T}^{0,x} \ \varepsilon \sigma(t,x)\\
    &\quad+ E[
       g(X_T^{t,x,0} + \varepsilon X_{1,T}^{t,x})H_{(i_1)}( X_{1,T}^{t,x},
      J_{t \to T}^{1,x}) 
    ] \ \varepsilon \sigma(t,x)\\
    &\quad+ \sum_{i_1,i_2} E[ g(X_T^{t,x,0} + \varepsilon X_{1,T}^{t,x}) H_{(i_2)}( X_{1,T}^{t,x}, H_{(i_1)}( X_{1,T}^{t,x}, X_{2,T}^{t,x} ) )      
    ]J_{t\to T}^{0,x} \ \varepsilon \sigma(t,x), 
  \end{align}
where the following relationship is taken into account: 
  \begin{align}
    H_{(i)}(X_T^{t,x,0} + \varepsilon X_{1,T}^{t,x}, 1)
    =H_{(i)}({X}_{1,T}^{t,x}, 1) / \varepsilon, \ \ \ i=1,\cdots,d. 
  \end{align} 
 Then, the similar calculation in (\ref{calc_ibp}) with (\ref{diff_computation})  gives the representation (\ref{cal V-1}). \ $\Box$

\section*{Acknowledgements}
This work is supported by JSPS KAKENHI (Grant Number 19K13736) and JST PRESTO (Grant Number JPMJPR2029), Japan.

\end{document}